\documentclass[12pt]{article}
\usepackage{graphicx}
\usepackage{amsmath}
\usepackage{amsfonts}
\usepackage{amssymb}
\newtheorem{theorem}{Theorem}

\newtheorem{lemma}[theorem]{Lemma}

\newtheorem{proposition}[theorem]{Proposition}
\newtheorem{remark}[theorem]{Remark}

\newenvironment{proof}[1][Proof]{\textbf{#1.} }{\ \rule{0.5em}{0.5em}}

\begin{document}

\title{TRANSVERSE RIEMANN-LORENTZ METRICS WITH TANGENT RADICAL }
\date{}
\author{E. Aguirre-Dab\'{a}n and J. Lafuente-L\'{o}pez\thanks{Departamento de
Geometr\'{i}a y Topolog\'{i}a, Facultad de Matem\'{a}ticas,
Universidad Complutense, 28040 Madrid, Spain. }} \maketitle
\begin{abstract}
Consider a smooth manifold with a smooth metric which changes
bilinear type from Riemann to Lorentz on a hypersurface $\Sigma$
with radical tangent to $\Sigma$. Two natural bilinear symmetric
forms appear there, and we use it to analyze the geometry of
$\Sigma$. We show the way in which these forms control the smooth
extensibility over $\Sigma$ of the covariant, sectional and Ricci
curvatures of the Levi-Civita connection outside
$\Sigma$.\smallskip
\newline \smallskip\textbf{Mathematics Subject Clasification (2000):} 53C50,
53B30, 53C15.\smallskip\newline \textbf{Key words: }Riemann-Lorentz metrics,
singular hypersurface, extensibility.
\end{abstract}

\section{Preliminaries}

Let $M$ be a $m$-dimensional connected manifold ($m>2)$ endowed with a smooth,
symmetric $(0,2)$-tensorfield $g$ which fails to have maximal rank on a (non
void) subset $\Sigma\subset M$. Thus, at each point $p\in\Sigma$, there exists
a nontrivial subspace (the \emph{radical}) $Rad_{p}\subset T_{p}M$, which is
orthogonal to the whole $T_{p}M$. We say that $(M,g)$ is a \emph{singular
space} (these spaces were analyzed for the first time in \cite{larsen}).
Moreover we say that $(M,g)$ is a \emph{transverse singular space} if, for any
local coordinate system $(x_{1},\ldots,x_{m})$, the function $\det
(g_{ab})_{a,b=1,\ldots,m}$ has non-zero differential at the points of $\Sigma$
(i.e. where $\det(g_{ab})$ vanishes). This implies at once: (1) the subset
$\Sigma$ is a smooth hypersurface in $M$, called the \emph{singular
hypersurface}; (2) at each point $p\in\Sigma$ the radical $Rad_{p}$ is one
dimensional; and (3) the signature of $g$ changes by $+1$ or $-1$ across
$\Sigma$ (see \cite{kossow1} for details). We say that $g$ has
\emph{transverse} (respectively \emph{tangent}) \emph{radical} if $Rad_{p}\cap
T_{p}\Sigma=\{0\}$ (respectively $Rad_{p}\subset T_{p}\Sigma$) for all
$p\in\Sigma$. There are several geometric and physical reasons to study such
spaces (see the Introduction to \cite{kossow1} for a detailed account) and
there are many articles devoted to the case with transverse radical (see
\cite{kossow4}, \cite{kossow1}, \cite{kossow2} and references therein).

In this article we analyze transverse singular spaces with \emph{tangent}
radical, and we focuse our attention on the case where \emph{the signature
changes from Riemann to Lorentz} across $\Sigma$.

In Section 2 we study the local geometry of the singular hypersurface $\Sigma
$, which is degenerate because the radical is assumed to be tangent. When the
degeneration of a hypersurface is due to its immersion in an overall
semiriemannian space, the Levi-Civita connection remains well-defined at the
points of the hypersurface (although it does not induce a connection on it);
in that case, and because the hypersurface has a one-dimensional normal vector
bundle which is moreover tangent to it, the geometry of the hypersurface (both
intrinsic and extrinsic) can be studied using the Levi-Civita connection,
looking for some (screen) distribution on the hypersurface complementary to
the normal bundle and, because the selected screen is non-canonical, focusing
the attention on those properties of the hypersurface which are
screen-independent (see e.g. \cite{bejancu}). But when the surrounding space
$(M,g)$ is singular at the points of the hypersurface (and independently of
whether the radical is transverse or tangent there), the Levi-Civita
connection fails to exist at such points; then, the suitable tool to analyze
the geometry of the hypersurface seems to be the canonically defined,
torsion-free, metric, dual connection on the whole $(M,g)$ (first pointed out
in \cite{kossow3}), which, in the case of one-dimensional radical, induces a
(conformally defined) symmetric $(0,2)$-tensorfield $II$ on the hypersurface.
This is what happens with our degenerate hypersurface $\Sigma$. We point out
that, because the normal vector bundle is now two-dimensional, there exists a
(local, determined up to a sign) \emph{canonical} smooth transverse
vectorfield $N$ along $\Sigma$ which is normal, unitary and $II$-isotropic.
This vectorfield $N$ allows us to construct, following the classical scheme, a
second fundamental form $\mathcal{H}$ on $\Sigma$, which in turn gives rise to
a \emph{canonical} screen distribution $S$ and also to a \emph{canonical}
vectorfield $R$ in the radical distribution. Thus, several nice canonical
structures arise in this case of tangent radical. All vectorfields tangent to
$\Sigma$ are uniquely decomposable in $S$- and $R$-components. With that
machinery, we analyze a natural family of torsion-free connections on $\Sigma
$, which we call \emph{admissible}; in case of $II$-flatness (i.e. when the
tensorfield $II$ vanishes), all such connections are metric and have the same
covariant curvatures.

In Section 3 we analyze the behaviour of some well-defined semiriemannian
objects (covariant derivatives, curvatures, ... on $M-\Sigma$), when we
approach the singular hypersurface $\Sigma$. We first point out that, by a
theorem in \cite{kossow1}, the transverse, $\ II$-isotropic vectorfield $N$
along $\Sigma$ has a \emph{canonical} (local) extension $\mathcal{N}$ to $M$
which is Levi-Civita geodesic outside $\Sigma$, and we use the flow of this
extension to (locally) extend to $M$ every vectorfield defined on $\Sigma$. We
apply these constructions to analyze whether or not the above mentioned
semiriemannian objects have good limits on $\Sigma$ and, when this is the
case, whether or not these limits only depend on the restriction to $\Sigma$
of the vectorfields we started with. We show how the second fundamental form
$\mathcal{H}$ and the symmetric $(0,2)$-tensorfield $II$ control these limit
properties. We compare our results with the corresponding ones \cite{kossow2}
in the case of transverse radical; the conclusion is that the tangent radical
case offers a wider variety of results and gives rise to some unavoidable
divergences (which are not present when the radical is transverse). In case of
$II$-flatness, we show a Gauss-type equation relating the curvature of the
admissible connections on $\Sigma$ with the good limit on $\Sigma$ of the
Levi-Civita curvature outside $\Sigma$.

It has been argued \cite{kossow5} that singular hypersurfaces with transverse
radical could provide a classical picture of the ''birth of time'' in general
relativity. Let us suppose that the radical happens to be tangent on some
(closed, with non void interior) region of the singular hypersurface $\Sigma$,
and that outside that region the radical recovers its generic behaviour and
becomes transverse. In view of the divergences occurring near the hypersurface
over the region with tangent radical, it would be interesting to analyze the
\emph{matching} of both regimes (e.g., across some compact submanifold of
$\Sigma$) and to look for traces (left by the tangent radical region) near the
hypersurface over the transverse radical region. This is beyond the scope of
the present article.

\subsection{Notations and conventions \label{notations}}

Vectorfields on $M$ are denoted by calligraphic letters $\mathcal{A}%
,\mathcal{B},\mathcal{C}, \ldots\in\frak{X}(M)$; we use $\mathcal{X}%
,\mathcal{Y},\mathcal{Z}, \ldots$ to denote vectorfields on $M$ tangent to
$\Sigma$. Vectorfields along $\Sigma$ are denoted by capital letters $A,B,C,
\ldots\in\frak{X}_{\Sigma}$; when they are tangent to $\Sigma$, we write
$X,Y,Z,\ldots\in\frak{X}(\Sigma)$. Note that $\frak{X}(\Sigma)$ is a
$C^{\infty}(\Sigma)$-submodule of $\frak{X}_{\Sigma}$. Given $\mathcal{A}
\in\frak{X}(M),$ we denote $\mathcal{A}^{o} \equiv\mathcal{A} \mid_{M -
\Sigma}\; \in\frak{X}(M-\Sigma)$ and $A \equiv\mathcal{A}\mid_{\Sigma}\;
\in\frak{X}_{\Sigma}$. In that case, we say that $\mathcal{A}$ is an
\emph{extension} of $A$.

\smallskip

Let us consider some function $\tau\in C^{\infty}(M)$ such that $\,\tau
\mid_{\Sigma}\;=0$ and $\,d\tau\mid_{\Sigma}\;\neq0$. We say that (locally,
around $\Sigma$) \emph{$\tau=0$ is an equation for $\Sigma$}. Given $f\in
C^{\infty}(M),$ it holds: $\;f\mid_{\Sigma}\;=0\;\Leftrightarrow\;f=k_{f}%
\tau\,,$ for some $k_{f}\in C^{\infty}(M)$. Denote $\,\tau^{o}\equiv\tau
\mid_{M-\Sigma}.$ When $f\mid_{\Sigma}\;=0,$ we shall write $\;{\tau^{o}}%
^{-1}f^{o}\cong0\;$ (the function $\;{\tau^{o}}^{-1}f^{o}\in C^{\infty
}(M-\Sigma)$ ''has no divergences when approaching $\Sigma$'' or ''has a good
limit on $\Sigma$'') and we shall say that \emph{$\;\tau^{-1}f$ is
well-defined} (as an element of $C^{\infty}(M)$).

\smallskip

Let $(\mathcal{E}_{1},\ldots,\mathcal{E}_{m})$ be a (local, $C^{\infty}$)
$\frak{X}(M)$-basis around some point of $\Sigma$. We say that $(\mathcal{E}%
_{1},\ldots,\mathcal{E}_{m})$ is an \emph{adapted frame} if $\mathcal{E}%
_{m}\mid_{\Sigma}\;$ is in the radical distribution and it holds
(orthonormality): $\;g_{ab}=\delta_{ab}(\pm(1-\delta_{am})+\delta_{am}%
\tau),\;$ for some $\tau\in C^{\infty}(M)\;(a,b=1,\ldots,m)$. Thus $\tau
\equiv\left\langle \mathcal{E}_{m},\mathcal{E}_{m}\right\rangle =0$ is an
equation for $\Sigma$. The existence of adapted frames around any point
$p\in\Sigma$ can be easily proved, e.g. starting with an orthonormal basis
$(e_{1},\ldots,e_{m})$ of $T_{p}M$ with $e_{m}\in Rad_{p}$, using a local
chart of $M$ around $p$ adapted to $\Sigma$, and applying a slight
modification of the Gram-Schmidt orthonormalization procedure. If $g$ has
transverse (respectively tangent) radical, all $\mathcal{E}_{i}^{\prime
}s\;(i=1,\ldots,m-1)$ can be chosen to be tangent (respectively, one of the
$\mathcal{E}_{i}^{\prime}s$ must be transverse) to $\Sigma$. \smallskip

We usually write $\,\left\langle \mathcal{A},\mathcal{B}\right\rangle ,$
instead of $g (\mathcal{A},\mathcal{B})$. We denote by $\nabla$ the
Levi-Civita connection outside $\Sigma$, and by $R$ its curvature. Connections
on $\Sigma$ are denoted by $D, \tilde{D},\dot{D},\ldots$, and their curvatures
by $R^{D}, R^{\tilde{D}}, R^{\dot{D}}, \ldots$.

\subsection{The Dual Connection \label{dual}}

Given a singular space $(M,g)$, there exists \cite{kossow3} a unique
torsion-free metric \emph{dual connection} on $M$, i.e. a unique map $\square:
\frak{X}(M) \times\frak{X} (M) \rightarrow\frak{X}^{\ast}(M)$ satisfying (for
all $f \in C^{\infty}(M)$ and $\mathcal{A},\mathcal{B},\mathcal{C} \in
\frak{X}(M)$):

(i) $\square_{f\mathcal{A}}\mathcal{B}= f \square_{\mathcal{A}} \mathcal{B}$

(ii) $\square_{\mathcal{A}} (f\mathcal{B}) = \mathcal{A}(f) \left\langle
\mathcal{B},.\right\rangle + f \square_{\mathcal{A}} \mathcal{B}$

(iii) $\square_{\mathcal{A}} \mathcal{B} (\mathcal{C}) - \square_{\mathcal{B}}
\mathcal{A} (\mathcal{C}) = \left\langle [\mathcal{A},\mathcal{B}%
],\mathcal{C}\right\rangle $ (torsion free)

(iv) $\square_{\mathcal{A}} \mathcal{B} (\mathcal{C}) + \square_{\mathcal{A}}
\mathcal{C} (\mathcal{B}) = \mathcal{A} \left\langle \mathcal{B}%
,\mathcal{C}\right\rangle $ (metric)

Moreover it holds:

(v) $\square_{\mathcal{A}^{o}}\mathcal{B}^{o} (\mathcal{C}^{o}) = \left\langle
\nabla_{\mathcal{A}^{o}} \mathcal{B}^{o},\mathcal{C}^{o} \right\rangle \;$
($\square$ is compatible with the Levi-Civita connection $\nabla$ on
$M-\Sigma$).

\smallskip

The dual connection can be alternatively defined by%

\begin{equation}
\left.
\begin{array}
[c]{l}%
2\square_{\mathcal{A}} \mathcal{B} (\mathcal{C}) := \mathcal{A} \left\langle
\mathcal{B},\mathcal{C} \right\rangle + \mathcal{B}\left\langle \mathcal{C}%
,\mathcal{A}\right\rangle -\mathcal{C}\left\langle \mathcal{A},\mathcal{B}%
\right\rangle +\\
\;\;\;\;\;\;\;\;\;\;\;\;\;\;\;\;\;\;\;\;\;\;\;\;\;\;\;\;\;\; + \left\langle
[\mathcal{A},\mathcal{B}],\mathcal{C}\right\rangle - \left\langle
[\mathcal{B},\mathcal{C}],\mathcal{A}\right\rangle + \left\langle
[\mathcal{C},\mathcal{A}],\mathcal{B}\right\rangle \;\;.
\end{array}
\right.  \label{square}%
\end{equation}

We realize that, $\forall\mathcal{A},\mathcal{B},\mathcal{C}\in\frak{X}%
(M),\;\square_{\mathcal{A}}\mathcal{B}(\mathcal{C})$ is $C^{\infty}(M)$-linear
in $\mathcal{A}$ and $\mathcal{C}$, thus $\square_{A}\mathcal{B}\in
\frak{X}_{\Sigma}^{\ast}$ is well-defined (we denote $A\equiv\mathcal{A}%
\mid_{\Sigma}$). This implies that:

(1) The dual connection has a good restriction $\square: \frak{X}(\Sigma)
\times\frak{X} (\Sigma) \rightarrow\frak{X}^{\ast}(\Sigma)$, which can also be
characterized as the unique torsion-free metric dual connection existing on
the singular hypersurface $\Sigma$.

(2) Given \emph{any} vectorfield $R\in\frak{X}_{\Sigma}$ in the radical
distribution, $\square_{A}\mathcal{B}(R)$ depends only on $A$ and
$B\equiv\mathcal{B}\mid_{\Sigma}$, thus $\square_{A}B(R)\in C^{\infty}%
(\Sigma)$ becomes well-defined and we obtain \cite{kossow3} a $C^{\infty
}(\Sigma)$-bilinear map $\,II_{R}:\frak{X}_{\Sigma}\times\frak{X}_{\Sigma
}\rightarrow C^{\infty}(\Sigma),\;(A,B)\mapsto\square_{A}B(R),$ which is
moreover symmetric because of (iii) above (see \cite{kossow2} for details). In
a similar way, given \emph{any} vectorfield $N\in\frak{X}_{\Sigma}$ orthogonal
to $\Sigma$, we obtain a $C^{\infty}(\Sigma)$-bilinear, symmetric map
$\;\mathcal{H}_{N}:\frak{X}(\Sigma)\times\frak{X}(\Sigma)\rightarrow
C^{\infty}(\Sigma),\;(X,Y)\mapsto\square_{X}Y(N).$ We shall come back to these
constructions later on.

\bigskip

\section{On the Geometry of the singular hypersurface}

From now on, we only consider transverse singular spaces $(M,g)$ with
\emph{tangent} radical.

\subsection{Fundamental forms}

At each point $p \in\Sigma$, the $g$-orthogonal subspace $T^{\perp}_{p}%
\Sigma\subset T_{p}M$ is a 2-plane and it holds: $T^{\perp}_{p} \Sigma\cap
T_{p}\Sigma= Rad_{p}$. Let $\,(\mathcal{E}_{1}, \ldots, \mathcal{E}_{m})$ be
an adapted frame around $p$, thus $E_{m}$ is in the radical distribution (we
denote $E_{a} \equiv\mathcal{E}_{a}\mid_{\Sigma}\, \in\frak{X}_{\Sigma
}\;,\;a=1,\ldots,m$), $\;\tau\equiv\left\langle \mathcal{E}_{m},\mathcal{E}%
_{m} \right\rangle = 0$ is an equation for $\Sigma$ and $E_{1}$ is transverse
to $\Sigma$. Therefore:%

\[
\left\{
\begin{array}
[c]{l}%
2II_{E_{m}} (E_{m}(p),E_{m}(p))= E_{m}(p) \left\langle \mathcal{E}%
_{m},\mathcal{E}_{m} \right\rangle = 0\\
2II_{E_{m}} (E_{1}(p),E_{m}(p))= E_{1}(p)\left\langle \mathcal{E}%
_{m},\mathcal{E}_{m} \right\rangle \neq0
\end{array}
\right.
\]

It follows that $II_{E_{m}}$ turns $T^{\perp}_{p}\Sigma$ into a Lorentz plane.
One of the two $\,II_{E_{m}}$-null directions at $p$ is determined by
$E_{m}(p)$. The other one cannot be a $\,g$-null direction, thus it determines
a unique (up to a sign) $g$-unitary vector in $T_{p}M$ $g$-orthogonal to
$\Sigma$. Moving from $p$ to the neighboring points in $\Sigma$, we obtain:

\begin{proposition}
Let $(M,g)$ be a transverse singular space with tangent radical and singular
hypersurface $\Sigma$. Then there exists locally a canonical smooth
vectorfield $N\in\frak{X}_{\Sigma}$ (determined up to a sign), which is
$g$-orthogonal to $\Sigma$, $II_{Rad}$-isotropic and $g$-unitary, i.e.
\[
\left\langle N,T\Sigma\right\rangle =0\;,\;\;\;II_{Rad}%
(N,N)=0\;\;\;\mathrm{and}\;\;\;\left\langle N,N\right\rangle =\pm1\;.\;\;\;
\]
If $g$ changes from Riemann to Lorentz, it holds: $\left\langle
N,N\right\rangle =1.$ We call $N$ the \emph{main normal} to $\Sigma
\;\;\;\;\;\;\;\;\;\;\endproof$
\end{proposition}

Over the singular hypersurface $\Sigma\;$ we have a \emph{first (degenerate)
fundamental form}, namely the restriction $g\mid_{\Sigma}$. Using now the main
normal $N$ to $\Sigma$, we are ready to construct (as in the classical theory)
a \emph{second fundamental form} $\mathcal{H}$ ($\equiv\mathcal{H}_{N}$,
already mentioned in \ref{dual}) \emph{on $\Sigma$}. We define, for
$X,Y\in\frak{X}(\Sigma)$:%

\[
\mathcal{H}(X,Y):=-\square_{X}N(Y)=\square_{X}Y(N)\;\;,
\]

\noindent where the second equality is because of property (iv) of the dual
connection. Moreover, property (iii) leads to the conclusion that
$\mathcal{H}$ is a symmetric $(0,2)$-tensor field over $\Sigma$. Note that, as
in the classical hypersurface theory, $\mathcal{H}$ is (locally) determined up
to a sign.

\smallskip

It turns out that, at each $p \in\Sigma$, the $\,II_{E_{m}}$-null direction
determined by $E_{m}(p)$ cannot be a $\mathcal{H}$-null direction, for it holds:%

\[
2\mathcal{H}(E_{m}(p),E_{m}(p)) = - 2II_{E_{m}}(E_{m}(p),N(p)) = -N(p)(\tau)
\neq0 \;\;;
\]

\noindent thus we can select a $\mathcal{H}$-unitary vectorfield $R\in
\frak{X}(\Sigma)$ in the radical distribution such that $\mathcal{H}%
(R,R)=\pm1$. \emph{Choosing $N$ such that $\mathcal{H}(R,R)=-1$}, we obtain

\begin{proposition}
Let $(M,g)$ be a transverse singular space with tangent radical and singular
hypersurface $\Sigma$. Then there exists locally a canonical smooth
vectorfield $R\in\frak{X}(\Sigma)$ (determined up to a sign), which is in the
radical distribution and $\mathcal{H}$-unitary, i.e.
\[
R(p)\in Rad_{p}\;,\;\mathrm{for\;all}\;p\in\Sigma\;,\;\;\;\mathrm{and}%
\;\;\;\mathcal{H}(R,R)=-1\;.\;\;\;
\]
We call $R$ the \emph{main radical vectorfield} on $\Sigma
\;\;\;\;\;\;\;\;\;\;\endproof$
\end{proposition}

\smallskip

The main radical vectorfield $R$ induces a canonical $C^{\infty}(\Sigma
)$-bilinear symmetric map%

\begin{equation}
II : \frak{X}_{\Sigma}\times\frak{X}_{\Sigma}\rightarrow C^{\infty}%
(\Sigma),\;(A,B) \mapsto\square_{A} B (R)\;\;\; \label{II}%
\end{equation}

\noindent(thus $II \equiv II_{R}$, already mentioned in \ref{dual}), whose
restriction to $\frak{X}(\Sigma) \times\frak{X}(\Sigma)$ yields \emph{another
symmetric $(0,2)$-tensorfield $II$ on $\Sigma$}. Note that it holds:%

\begin{equation}
\left\{
\begin{array}
[c]{l}%
II(N,N) = 0\\
II(N,X) = - \mathcal{H}(X,R)\;\;,\;\;\mathrm{for\;all}\;X \in\frak{X}%
(\Sigma)\\
II(A,R) = 0 \;\;\Leftrightarrow\;\;A \in\frak{X}(\Sigma) \;\;
\end{array}
\right.  \;\;\; \label{HII}%
\end{equation}

\subsection{The canonical Screen Distribution}

A \emph{screen distribution} is a distribution on the singular hypersurface
$\Sigma$ which yields, at each $p\in\Sigma$, a hyperplane of $T_{p}\Sigma$
transversal to $Rad_{p}$. We now define the \emph{canonical screen
distribution} by choosing (at each $p\in\Sigma$)%

\begin{equation}
S_{p}:=\left\{  v\in T_{p}\Sigma\;:\;\mathcal{H}(v,R(p))=0\right\}
\;\;,\;\;\; \label{S}%
\end{equation}
where $R$ is the main radical vectorfield. We shall denote by $S$ either the
set $\left\{  S_{p}:p\in\Sigma\right\}  \;$ or the corresponding vector bundle
(the \emph{screen bundle}) $\;S\rightarrow\Sigma$. We denote by $\Gamma(S)$
the $C^{\infty}(\Sigma)$-module of sections of $S$, which is a submodule of
$\frak{X}(\Sigma)$.

\emph{From now on, we only consider the case where the signature of $g$
changes from Riemann to Lorentz upon crossing $\Sigma$}\footnote{Instead of
that hypothesis, we could also make the (weaker) requirement that the metric
$g$ is non-degenerate over any slice of the screen distribution. Then the
screen bundle would become a semiriemannian vector bundle. All results in this
article would remain essentially valid.}. This means that $g$ is semi-defined
on $\Sigma$ and the screen bundle becomes a riemannian vector bundle.

The proof of the following proposition is straightforward.

\smallskip

\begin{proposition}
Let $(M,g)$ be a transverse Riemann-Lorentz space with tangent radical,
singular hypersurface $\Sigma$ and canonical screen distribution $S$. Then
there exist, over each euclidean fibre $S_{p}$ of $S$, two selfadjoint
endomorphisms $\mathcal{H}_{p}^{S}$ and $II_{p}^{S}$ such that (for all
$v,w\in S_{p}$):%

\[
\left\{
\begin{array}
[c]{l}%
\left\langle \mathcal{H}_{p}^{S}(v),w\right\rangle :=\mathcal{H}(v,w)\\
\left\langle II_{p}^{S}(v),w\right\rangle :=II(v,w)
\end{array}
\right.  \;\;;
\]

\noindent these determine vector bundle endomorphisms $\mathcal{H}^{S}$ and
$II^{S}$ of $S$, called the \emph{Weingarten screen maps}. The eigenvalues of
these screen maps are the \emph{principal curvatures}, and the correspoding
eigenvectors define the \emph{principal directions} of $\mathcal{H}^{S}$ or
$II^{S}$, respectively $\;\;\;\;\;\;\;\;\;\;\endproof$
\end{proposition}

\bigskip

The natural definition of \emph{$II$-flatness} (respectively,
\emph{$\mathcal{H}$-flatness}) of $\Sigma$ is that $II^{S}$ (respectively,
$\mathcal{H}^{S}$) is identically zero; or, in other words:%

\begin{equation}
II(V,W)=0\;\;\;\;\;(\mathrm{respectively,}\;\;\mathcal{H}%
(V,W)=0)\;,\;\mathrm{for\;all}\;V,W\in\Gamma(S)\;\;. \label{IIflat}%
\end{equation}

Because of (\ref{HII}), $II$-flatness is equivalent to $II(X,Y)=0$, for all
$X,Y \in\frak{X}(\Sigma)$. And, because of (\ref{S}), $\mathcal{H}$-flatness
is equivalent to $\mathcal{H}(V,X)=0$, for all $V \in\Gamma(S)$ and $X
\in\frak{X}(\Sigma)$. \bigskip

\begin{remark}
\label{rem:IIflat} A few words about the definition of $II$-flatness. As we
mentioned in \ref{dual}, given any transverse singular space $(M,g)$ with
singular hypersurface $\Sigma$, and for any nowhere zero vectorfield
$R\in\frak{X}(\Sigma)$ in the radical distribution, a well-defined $C^{\infty
}(\Sigma)$-bilinear, symmetric map $\,II:\frak{X}_{\Sigma}\times
\frak{X}_{\Sigma}\rightarrow C^{\infty}(\Sigma),(A,B)\mapsto\square_{A}B(R)\,$
arises. In this general context, $\Sigma$ was defined \cite{kossow2} to be
$II$-flat if it holds:%

\[
II(A,B)=0\;,\;for\;all\;\;B\in\frak{X}_{\Sigma}\;\;\;\Leftrightarrow
\;\;\;A\in\frak{X}(\Sigma)\;\;.
\]

Now it is straightforward to see (because $\;II(A,R)=0\;\Leftrightarrow
\;A\in\frak{X}(\Sigma)$) that this requirement is equivalent to the following
two conditions:%

\[
\left\{
\begin{array}
[c]{l}%
(i)\;II(X,Y)=0\;\;,\;\;for\;all\;\;X,Y\in\frak{X}(\Sigma)\\
(ii)\;the\;\;radical\;\;is\;\;transverse\;\;.
\end{array}
\right.
\]

Because $(i)$ and $(ii)$ are in fact independent conditions, our definition of
$II$-flatness (which is just condition $(i)$) turns out to be the natural one
for general case (and coincides with the definition given in \cite{kossow2}
for the case with transverse radical) $\;\;\;\;\;\endproof$
\end{remark}

A vectorfield $A \in\frak{X}_{\Sigma}$ can be decomposed in normal-, screen-
and radical-components, as follows%

\[
A=\nu(A)N+A^{S}+\rho(A)R\;\;,\;\;\;
\]

\noindent where $\nu(A):=\left\langle A,N\right\rangle $ and $\rho
(A):=-\mathcal{H}(A-\nu(A)N,R)$. Thus $\rho\in\frak{X}_{\Sigma}^{\ast}$ is
completely determined by the (equally denoted) 1-form $\rho=-\mathcal{H}%
(\,.\,,R)\in\frak{X}^{\ast}(\Sigma)$. Because $\;d\rho(V,W)=-\rho([V,W]),$ for
all $V,W\in\Gamma(S)$, it follows: $\;d\rho=0\;\Rightarrow\;$ the distribution
$\,S$ is integrable.

\subsection{Admissible Connections}

We are going to describe some natural connections on $\Sigma$. All these
connections arise without any explicite reference to the Levi-Civita
connection outside $\Sigma$. We first analyze a canonical \emph{screen
connection-operator}, which yields a natural metric connection on the
riemannian screen vector bundle $\;S\rightarrow\Sigma$. Requiring
compatibility with that operator and zero torsion leads to a family of
connections which we call \emph{admissible}. These connections are not
necessarily metric. It turns out that, if there exists a torsion-free metric
connection on $\Sigma$, then: (1) it is necessarily admissible; (2) all
admisible connections are metric; and (3) $\Sigma$ is $II$-flat. We finally
prove that $II$-flatness is also the necessary and sufficient condition for
all admissible connections induce the same covariant curvature tensor. In
Section \ref{near}, starting with the Levi-Civita connection outside $\Sigma$
(and provided that $\Sigma$ is $II$-flat), we shall induce on $\Sigma$ (like
in the classical theory of semiriemannian hypersurfaces) the so called
\emph{tangential connection}, which turns out to be admissible.

\bigskip

Because $g\mid_{S}$ does not degenerate, it is natural to introduce the
\emph{screen connection-operator} as the map $\;D^{S}:\frak{X}(\Sigma
)\times\frak{X}_{\Sigma}\rightarrow\Gamma(S),(X,A)\mapsto D_{X}^{S}A$, defined
by
\[
(D_{X}^{S}A)(p):=D_{X(p)}^{S}A\;\;\;,\;\;\;\mathrm{for\;all}\;p\in\Sigma\;\;,
\]

\noindent where $\,D^{S}_{X(p)} A \in S_{p} \subset T_{p}\Sigma\;$ is the
\emph{unique} vector satisfying%

\[
\square_{X(p)} A (v) = \left\langle D^{S}_{X(p)} A, v \right\rangle
\;\;\;,\;\;\;\mathrm{for\;all}\;v\in S_{p}\;\;.
\]

The screen connection-operator $D^{S}$ has the following properties (for all
$f \in C^{\infty}(\Sigma),\;X \in\frak{X}(\Sigma),\;A,B \in\frak{X}_{\Sigma}
\;\;\mathrm{and}\;\; V \in\Gamma(S)$):

(i) $II^{S}(V)=-D^{S}_{V} R$

(ii) $\mathcal{H}^{S}(V)=-D^{S}_{V} N$

(iii) $D^{S}_{fX}A = f D^{S}_{X} A \;\;\;\mathrm{and}\;\;\;D^{S}_{X+Y} A =
D^{S}_{X} A + D^{S}_{Y} A$

(iv) $D^{S}_{X} (A+B) = D^{S}_{X} A + D^{S}_{X} B$

(v) $D^{S}_{X} (fA) = X(f)A^{S} + f D^{S}_{X} A$; in particular, $D^{S}_{X}
(fR) = f D^{S}_{X} R$

(vi) $\left\langle D^{S}_{X} V, W \right\rangle + \left\langle V , D^{S}_{X} W
\right\rangle = X \left\langle V,W \right\rangle \;\;\;\;$ (metric)

(vii) $\left\langle D^{S}_{X} A,V \right\rangle = \square_{X} A (V)
\;\;\;\;(\mathrm{compatible\;with}\;\square)\;\;.$

\smallskip

Properties (iii)-(vi) show that $\;D^{S} : \frak{X}(\Sigma) \times\Gamma(S)
\rightarrow\Gamma(S)\;$ gives a metric connection on the screen riemannian
vector bundle $S \rightarrow\Sigma$.

However, property (v) shows that the restriction $\,D^{S} : \frak{X}(\Sigma)
\times\frak{X}(\Sigma) \rightarrow\Gamma(S)\,$ \emph{does not} give a
connection on $\Sigma$.

\bigskip

When looking for interesting connections $D$ on $\Sigma$, it is natural to
impose two requirements: (1) $D$ should satisfy (for all $X,Y\in
\frak{X}(\Sigma)$): $\;(D_{X}Y)^{S}=D_{X}^{S}Y\,,$ or equivalently:
$\left\langle D_{X}Y,V\right\rangle =\square_{X}Y(V),$ for all $V\in\Gamma
(S)$; and (2) $D$ should be torsion-free. We shall call such connections
\emph{admissible}.

\smallskip

The most obvious connection on $\Sigma$ satisfying the first condition is the
following one (we denote it by $\tilde{D}$ and we call it the \emph{main
connection})%

\begin{equation}
\tilde{D}:\;\frak{X}(\Sigma) \times\frak{X}(\Sigma) \rightarrow\frak{X}%
(\Sigma)\;,\;(X,Y) \; \mapsto\; \tilde{D}_{X} Y := D^{S}_{X} Y + X(\rho(Y))
R\;,\;\; \label{main}%
\end{equation}

\noindent whose properties are analyzed in the next

\begin{proposition}
\label{prop:main} Let $(M,g)$ be a transverse Riemann-Lorentz space with
tangent radical and singular hypersurface $\Sigma$. Then, the main connection
$\tilde{D}$ on $\Sigma$:

(a) has torsion $\widetilde{Tor}=R\otimes d\rho\;$.

(b) satisfies (for all $X,Y,Z\in\frak{X}(\Sigma)$): $\;\left\langle \tilde
{D}_{X}Y,Z\right\rangle +\left\langle Y,\tilde{D}_{X}Z\right\rangle
=X\left\langle Y,Z\right\rangle \,$ (i.e. is metric) if and only if $\Sigma$
is $II$-flat.
\end{proposition}

\begin{proof}
In what follows, $S$ is the canonical screen distribution.

(a) Let be $X,Y\in\frak{X}(\Sigma)$. By Property (vii) of $D^{S}$, one
immediately sees that:$\;(\widetilde{Tor}(X,Y))^{S}=0.$ Therefore,
$\;\widetilde{Tor}(X,Y)=\rho(\widetilde{Tor}(X,Y))R=d\rho(X,Y)R$.

(b) We have:%

\[
\left.
\begin{array}
[c]{l}%
\left\langle \tilde{D}_{X}Y,Z\right\rangle =\left\langle D_{X}^{S}%
Y,Z^{S}\right\rangle =\left\langle D_{X}^{S}Y^{S},Z^{S}\right\rangle
+\rho(Y)\left\langle D_{X}^{S}R,Z^{S}\right\rangle =\\
\;\;\;\;\;\;\;\;\;\;=\;\;(\mathrm{by\;property\;(vii)\;of}\;D^{S}%
)\;\;\left\langle D_{X}^{S}Y^{S},Z^{S}\right\rangle -\rho(Y)II(X,Z)\;;
\end{array}
\right.
\]

\noindent thus, by property (vi) of $D^{S}$, we obtain:%

\[
\left\langle \tilde{D}_{X}Y,Z\right\rangle +\left\langle Y,\tilde{D}%
_{X}Z\right\rangle =X\left\langle Y,Z\right\rangle -II(X,\rho(Y)Z+\rho(Z)Y)\;.
\]

Now ($\Leftarrow$) is trivial. Let us prove ($\Rightarrow$): If $\tilde{D}$ is
metric, last formula yields: $\;\rho(X)II(X,X)=0,\;$ for all $X\in
\frak{X}(\Sigma)$. Because $T_{p}\Sigma-S_{p}$ is dense in $T_{p}\Sigma$ (for
all $p\in\Sigma$), it follows that $\,II(X,X)=0,\;$ for all $X\in
\frak{X}(\Sigma)$. Being $II$ symmetric, this implies that $\Sigma$ is
$II$-flat $\;\;\;\;$
\end{proof}

\bigskip

Thus, unless $\;d\rho=0,$ the main connection $\tilde{D}$ is not admissible.
However, it is straightforward to check that the connection%

\begin{equation}
\dot{D}:\frak{X}(\Sigma) \times\frak{X}(\Sigma) \rightarrow\frak{X}%
(\Sigma),\;(X,Y) \mapsto\; \dot{D}_{X} Y := \tilde{D}_{X} Y - \frac{1}{2} d
\rho(X,Y) R\; \label{main admissible}%
\end{equation}

\noindent is always admissible (we call it the \emph{main admissible
connection}).

Admissible connections have the following properties:

\begin{theorem}
\label{theor:admissible} Let $(M,g)$ be a transverse Riemann-Lorentz space
with tangent radical and singular hypersurface $\Sigma$. Then:

(a) Any admissible connection $D$ on $\Sigma$ satisfies (for all
$X,Y\in\frak{X}(\Sigma)$)%

\[
D_{X}Y=\dot{D}_{X}Y+\sigma(X,Y)R\;\;,
\]

\noindent where $\dot{D}$ is the main admissible connection and $\sigma$ is
some symmetric $(0,2)$-tensorfield on $\Sigma$.

(b) If there exists a torsion-free metric connection on $\Sigma$, then: (1) it
is admissible; (2) all admissible connections on $\Sigma$ are metric; and (3)
$\Sigma$ is $II$-flat.

(c) If $\Sigma$ is $II$-flat, all admissible connections have the same
covariant curvature.
\end{theorem}

\begin{proof}
In what follows, $X,Y,Z$ are arbitrary in $\frak{X}(\Sigma)$ and $V$ is
arbitrary in $\Gamma(S)$, where $S$ is the canonical screen distribution.

(a) As is well known, any torsion-free connection $D$ on $\Sigma$ must
satisfy: $\;D_{X}Y=\dot{D}_{X}Y+\varphi(X,Y),\;$ where $\varphi$ is some
symmetric $(1,2)$-tensorfield on $\Sigma$. Being both $D$ and $\dot{D}$
admissible, it must hold: $\;\left\langle \varphi(X,Y),V\right\rangle =0,\;$
and the result follows.

(b) (1) If $D$ is a torsion-free metric connection on $\Sigma$, the induced
dual connection $\square^{D}$, defined by: $\;\square_{X}^{D}%
Y(Z):=\left\langle D_{X}Y,Z\right\rangle ,\;$ becomes torsion-free and metric;
thus, by uniqueness (see \ref{dual}), it must hold: $\square^{D}=\square$, and
$D$ becomes admissible. (2) If $D$ is an admissible connection on $\Sigma$, it
follows from Part (a) that: $\;D_{X}Y=\tilde{D}_{X}Y+(\sigma-\frac{1}{2}%
d\rho)(X,Y)R,\;$ for some symmetric $(0,2)$-tensorfield $\sigma$ on $\Sigma$.
But this implies: $D$ is metric if and only if $\tilde{D}$ is metric. (3)
follows from (2) and from Part (b) of Proposition \ref{prop:main}.

(c) Let $D$ be any admissible connection on $\Sigma$ and let us consider its
covariant curvature, defined (for all $X,Y,Z,T\in\frak{X}(\Sigma)$) by:%

\[
\left\langle R^{D}(X,Y)Z,T\right\rangle :=\left\langle D_{X}(D_{Y}%
Z)-D_{Y}(D_{X}Z)-D_{[X,Y]}Z,T\right\rangle \;.
\]

Using Part (a), we compute:%

\[
\left.
\begin{array}
[c]{c}%
\left\langle D_{X}(D_{Y}Z),T\right\rangle =\left\langle D_{X}(\dot{D}%
_{Y}Z+\sigma(Y,Z)R),T\right\rangle =\left\langle \dot{D}_{X}(\dot{D}%
_{Y}Z)+\sigma(Y,Z)\dot{D}_{X}R,T\right\rangle =\\
\\
=\left\langle \dot{D}_{X}(\dot{D}_{Y}Z),T\right\rangle +\sigma
(Y,Z)\left\langle D_{X}^{S}R,T\right\rangle =\left\langle \dot{D}_{X}(\dot
{D}_{Y}Z),T\right\rangle -\sigma(Y,Z)II(X,T)\;\;,
\end{array}
\right.
\]

\noindent where last equality is because: $\;\left\langle D_{X}^{S}%
R,T\right\rangle =\left\langle D_{X}^{S}R,T^{S}\right\rangle =\square
_{X}R(T^{S})=\square_{X}R(T)=-II(X,T)$. Therefore we obtain:%

\[
\left\langle R^{D}(X,Y)Z,T\right\rangle =\left\langle R^{\dot{D}%
}(X,Y)Z,T\right\rangle -\det\left(
\begin{array}
[c]{c}%
\sigma(Y,Z)\;\;II(Y,Z)\\
\sigma(X,T)\;\;II(X,T)
\end{array}
\right)  \;\;,
\]

\noindent and the result follows $\;\;\;\;\;$
\end{proof}

\section{Near the singular hypersurface \label{near}}

We analyze in this Section the behaviour of some well-defined Levi-Civita
objects (covariant derivatives, curvatures, ... on $M-\Sigma$) when we
approach the singular hypersurface $\Sigma$. Thus we can replace (if
necessary) the whole $M$ by a (small enough) neighborhood of $\Sigma$ in $M$.

Typically, we start with vectorfields $\mathcal{A},\mathcal{B},\ldots
\in\frak{X}(M)$, construct some semiriemannian object, say $\,\bigcirc
(\mathcal{A}^{o},\mathcal{B}^{o},\ldots)$, on $M-\Sigma$, and ask under what
circumstances this object: (1) has a good limit on $\Sigma$, i.e.
$\bigcirc(\mathcal{A}^{o},\mathcal{B}^{o},\ldots)\cong0$ (in that case we say
that ''$\bigcirc(\mathcal{A},\mathcal{B},\ldots)$ is well-defined''); and (2)
the restriction $\bigcirc(\mathcal{A},\mathcal{B},\dots)\mid_{\Sigma}\;$ only
depends on $A\equiv\mathcal{A}\mid_{\Sigma},\ldots\in\frak{X}_{\Sigma}$ (in
that case we say that ''$\bigcirc(A,B,\ldots)$ is well-defined'').

When dealing with two such objects $\,\bigcirc_{1}$ and $\,\bigcirc_{2}$, we
usually shall write $\,\bigcirc_{1}(\mathcal{A}^{o},\mathcal{B}^{o},..)
\cong\bigcirc_{2}(\mathcal{A}^{o},\mathcal{B}^{o},..)$ to mean $\,\bigcirc
_{1}(\mathcal{A}^{o},\mathcal{B}^{o},..) -\bigcirc_{2}(\mathcal{A}%
^{o},\mathcal{B}^{o},..) \cong0.$

\subsection{Vectorfield extensions}

Let $\nabla$ be the Levi-Civita connection outside $\Sigma$ and let us
consider the main normal $N \in\frak{X}_{\Sigma}$, thus $\left\langle N,N
\right\rangle = 1$. Because $II(N,N)= 0,$ it follows from Theorem 1 in
\cite{kossow1} that there exists a \emph{unique} (we call it \emph{canonical})
local extension $\mathcal{N} \in\frak{X}(M)$ of $N$ which is $\nabla$-geodesic
outside $\Sigma$. By continuity, it follows that: $\left\langle \mathcal{N}%
,\mathcal{N} \right\rangle = 1$.

\bigskip

Given $A \in\frak{X}_{\Sigma}$, there exists a \emph{unique} (we call it
\emph{canonical}) local extension $\mathcal{A} \in\frak{X}(M)$ such that%

\[
[\mathcal{N},\mathcal{A}]=0
\]

\noindent($\mathcal{A}$ is generated from $A$ by the flow of $\mathcal{N}$),
and for that extension it holds:%

\[
\mathcal{N}^{o} \left\langle \mathcal{N}^{o}, \mathcal{A}^{o} \right\rangle =
\left\langle \mathcal{N}^{o}, \nabla_{\mathcal{N}^{o}} \mathcal{A}^{o}
\right\rangle = \left\langle \mathcal{N}^{o}, \nabla_{\mathcal{A}^{o}}
\mathcal{N}^{o} \right\rangle = 0 \;,\;\Rightarrow\;\mathcal{N} \left\langle
\mathcal{N},\mathcal{A} \right\rangle = 0 \;;
\]

\noindent in particular, if $X \in\frak{X}(\Sigma)\;(\Rightarrow\;\left\langle
N,X \right\rangle = 0),$ the canonical extension $\mathcal{X}\in\frak{X} (M)$ satisfies:%

\begin{equation}
\;\left\langle \mathcal{N},\mathcal{X}\right\rangle =0\;\;\;. \label{NX}%
\end{equation}

\bigskip

Let $\mathcal{R} \in\frak{X}(M)$ be the canonical extension of the main
radical vectorfield $R$. In what follows, we shall denote $\;\tau
\equiv\left\langle \mathcal{R},\mathcal{R} \right\rangle $. We obtain from
(\ref{square}):%

\begin{equation}
N (\tau) = 2 \square_{R} N (R)=:-2\mathcal{H}(R,R) = 2\;;\;\; \label{Ntau}%
\end{equation}

\noindent thus $\;\tau=0\;$ is an equation for $\Sigma$ and $\;\tau
^{-1}(\mathcal{N}(\tau)-2)\,$ is well-defined (by the way, (\ref{Ntau}) is
still valid if $\mathcal{N},\mathcal{R}$ are arbitrary extensions of $N,R$).

\bigskip

Let be $X\in\frak{X}(\Sigma)$ and let $\mathcal{X}\in\frak{X}(M)$ be its
\emph{canonical} extension. It follows from $\,X(\tau)=0$ that $\;\tau
^{-1}\mathcal{X}(\tau)$ is well-defined and it holds:%

\[
\mathcal{N}(\mathcal{X}(\tau)) = \mathcal{X}(\mathcal{N}(\tau)) \equiv
\mathcal{X}([\tau^{-1}(\mathcal{N}(\tau) -2)]\;\tau) \;,\;\;\Rightarrow
\;\;N(\mathcal{X}(\tau)) = 0\;\;;
\]

\noindent thus $\;\tau^{-1}\mathcal{N}(\tau)$ is well-defined and it holds:%

\[
\left.
\begin{array}
[c]{c}%
(\tau^{-1}\mathcal{N}(\mathcal{X}(\tau))) \;\tau\equiv\mathcal{N}%
(\mathcal{X}(\tau)) \equiv\mathcal{N}([\tau^{-1}\mathcal{X}(\tau)]\;\tau) =\\
\\
= \mathcal{N}(\tau^{-1}\mathcal{X}(\tau))\tau+ (\tau^{-1}\mathcal{X}(\tau)) (2
+ [\tau^{-1}(\mathcal{N}(\tau) -2]) \;\tau)\;\;,\;\;\Rightarrow\;\;(\tau
^{-1}\mathcal{X}(\tau))\mid_{\Sigma}\;= 0 \;\;;
\end{array}
\right.
\]

\noindent therefore%

\begin{equation}
\tau^{-2}\mathcal{X}(\tau)\;\;\mathrm{is\;well-defined}\;\;.\;\; \label{Xtau}%
\end{equation}

\bigskip

Let be $X \in\frak{X}(\Sigma)$ and let $\mathcal{X} \in\frak{X}(M)$ be
\emph{any} extension. A direct computation from (\ref{square}) leads to (we
use (\ref{Ntau})):%

\[
\left.
\begin{array}
[c]{c}%
\mathcal{H}(X,R) := \square_{\mathcal{X}} \mathcal{R} (\mathcal{N})
\mid_{\Sigma}\; = -\frac{1}{2} (\mathcal{N} \left\langle \mathcal{X}%
,\mathcal{R} \right\rangle )\mid_{\Sigma}\;\equiv-\frac{1}{2} N ((\tau^{-1}
\left\langle \mathcal{X},\mathcal{R} \right\rangle )\;\tau) =\\
\\
= - \frac{1}{2} (\tau^{-1} \left\langle \mathcal{X},\mathcal{R} \right\rangle
)\mid_{\Sigma}\;N(\tau) = - (\tau^{-1} \left\langle \mathcal{X},\mathcal{R}
\right\rangle )\mid_{\Sigma}\;\;;
\end{array}
\right.
\]

\noindent therefore, given $A \in\frak{X}_{\Sigma}$, and for \emph{any}
extension $\mathcal{A} \in\frak{X}(M)$ of $A$, it holds:%

\begin{equation}
\rho(A):= - \mathcal{H}(A-\nu(A)N,R) = (\tau^{-1} \left\langle \mathcal{A}%
,\mathcal{R} \right\rangle ) \mid_{\Sigma}\;\;.\;\; \label{ro}%
\end{equation}

\bigskip

\emph{In what follows, we shall use some fixed ($C^{\infty}$, local around
some point of $\Sigma$) adapted frame $(\mathcal{E}_{1},\ldots,\mathcal{E}%
_{m}=\mathcal{R})$, with $\;E_{1}=N,E_{2},\ldots,E_{m-1}\in\Gamma(S)$ and
$\mathcal{R}$ the canonical extension of $R$}.

\begin{remark}
In this paper, we do not use coordinates. Nevertheless, it is interesting to
point out that we can construct, around each point of $\Sigma$, a coordinate
system $\,(x^{1},\ldots,x^{m})$ of $M$ such that (we call these coordinates
\emph{adapted}): (1) $\partial_{x^{1}}=\mathcal{N}$ and $\partial_{x^{m}%
}=\mathcal{R}$ (the canonical extensions of $N$ and $R$); (2) $x^{1}=0$ is an
equation for $\Sigma$; and (3) it holds:%

\[
(g_{ab})=\left(
\begin{array}
[c]{ccc}%
1 & 0 & 0\\
0 & (g_{\lambda\mu}) & x^{1}g_{\lambda}\\
0 & x^{1}g_{\lambda} & x^{1}g_{m}%
\end{array}
\right)  \;\;(a,b=1,\ldots,m\,;\,\lambda,\mu=2,\ldots,m-1),
\]

\noindent for some $\;g_{i}\in C^{\infty}(M)\;\;(i=2,\ldots,m)\;$ with
$\;g_{m}(0,x^{2},\ldots,x^{m})=2$.

\smallskip

We now outline the construction: Around $\;p\in\Sigma$, we choose coordinates
$\,(x^{2},\ldots,x^{m})$ in $\Sigma$ such that $\partial_{x^{m}}=R$. Using the
flow of $\mathcal{N}$, it is straightforward to construct a coordinate system
$\,(x^{1},\ldots,x^{m})$ in $M$ such that $\partial_{x^{1}}=\mathcal{N}$ (thus
$g_{11}=1$) and $\partial_{x^{i}}$ is the canonical extension of the (equally
denoted) vectorfield $\partial_{x^{i}}\in\frak{X}(\Sigma)\;(i=2,\ldots,m)$; in
particular, $\partial_{x^{m}}=\mathcal{R}$. Obviously, $x^{1}=0$ is an
equation for $\Sigma$. It follows from (\ref{NX}): $\;g_{1i}=0\;(i=2,\ldots
,m)$. On the other hand, $\;g_{im}\mid_{\Sigma}\;=0$ implies: $\;g_{im}%
=x^{1}g_{i},$ for some $\;g_{i}\in C^{\infty}(M)\;(i=2,\ldots,m)$. And
finally, writing $\;\tau\equiv g_{mm},$ it follows from (\ref{Ntau}):
$\;g_{m}\mid_{\Sigma}\;=2$.

Note that, in adapted coordinates, a vectorfield $\;\mathcal{A}\equiv\sum
f_{a}\partial_{x^{a}}\in\frak{X}(M)$ is the canonical extension of
$\;A\equiv\sum f_{a}\mid_{\Sigma}\partial_{x^{a}}\in\frak{X}_{\Sigma}$ if and
only if all $f_{a}^{\prime}s\;(a=1,\ldots,m)$ do not depend on $x^{1}.$

\smallskip

Calling $\;\Gamma_{cab}\equiv\square_{\partial_{x^{a}}}\partial_{x^{b}%
}(\partial_{x^{c}})$, the first-class Christoffel symbols of $\;\square$ in
that coordinates, and using (\ref{square}), it is straigthforward to see that
the components of the second fundamental form $\,\mathcal{H}\,$ and of the
tensorfield $\,II\,$ on $\Sigma$ are given by%

\[
\left\{
\begin{array}
[c]{c}%
\mathcal{H}_{ij}=\Gamma_{1ij}\mid_{\Sigma}\text{ }=-\frac{1}{2}{\frac{\partial
g_{ij}}{\partial x^{1}}}\mid_{x^{1}=0}\\
\\
II_{ij}=\Gamma_{mij}\mid_{\Sigma}\text{ }=-\frac{1}{2}{\frac{\partial g_{ij}%
}{\partial x^{m}}}\mid_{x^{1}=0}%
\end{array}
\right.  \;\;(i,j=2,\ldots,m)\;\;.
\]

\smallskip

The above construction applies to all local examples of transverse
Riemann-Lorentz metrics with tangent radical $\;\;\;\;\;\endproof$
\end{remark}

\subsection{Covariant derivatives}

Let be $\mathcal{A},\mathcal{B} \in\frak{X}(M)$. On $M-\Sigma\;$ we have:%

\[
\nabla_{\mathcal{A}^{o}} \mathcal{B}^{o} = \Sigma^{m-1}_{i=1} \square
_{\mathcal{A}^{o}} \mathcal{B}^{o} (\mathcal{E}^{o}_{i}) \mathcal{E}^{o}_{i} +
{\tau^{o}}^{-1} \square_{\mathcal{A}^{o}} \mathcal{B}^{o} (\mathcal{R}^{o})
\mathcal{R}^{o}\;.
\]

It thus follows that the vectorfield $\nabla_{\mathcal{A}^{o}} \mathcal{B}%
^{o}$ has a good limit on $\Sigma\;$ if and only if: $\tau^{-1} \square
_{\mathcal{A}} \mathcal{B} (\mathcal{R}) \in C^{\infty}(M)\,$ or, using
(\ref{II}), if and only if: $II(A,B)=0$.

Once $\nabla_{\mathcal{A}^{o}}\mathcal{B}^{o}\cong0$, we may define:%

\begin{equation}
\nabla_{\mathcal{A}} \mathcal{B} := \Sigma^{m-1}_{i=1} \square_{\mathcal{A}}
\mathcal{B} (\mathcal{E}_{i}) \mathcal{E}_{i} + \tau^{-1} \square
_{\mathcal{A}} \mathcal{B} (\mathcal{R}) \;\mathcal{R} \;\;,\;\;\;
\label{covar.deriv.}%
\end{equation}

\smallskip

\noindent which obviously satisfies%

\[
\left\langle \nabla_{\mathcal{A}} \mathcal{B},\mathcal{C} \right\rangle =
\square_{\mathcal{A}}\mathcal{B} (\mathcal{C})\;\;,\;\;\mathrm{for\;all}%
\;\;\mathcal{C} \in\frak{X}(M)\;\;.
\]

It turns out that, in general, the restriction $\nabla_{\mathcal{A}}
\mathcal{B} \mid_{\Sigma}\;$ depends, not only on the restrictions $A,B
\in\frak{X}_{\Sigma}$, but also on the original vectorfields $\mathcal{A}%
,\mathcal{B}$. Indeed, starting with $\mathcal{A}^{\prime}=\mathcal{A} +
\tau\bar{\mathcal{A}},\;\mathcal{B}^{\prime}=\mathcal{B} + \tau\bar
{\mathcal{B}}$ (for some $\bar{\mathcal{A}},\bar{\mathcal{B}} \in\frak{X}%
(M)$), it is straightforward to see that it holds (we denote $\bar{A}
\equiv\bar{\mathcal{A}}\mid_{\Sigma}\;,\bar{B} \equiv\bar{\mathcal{B}}%
\mid_{\Sigma}\;$):%

\[
(\tau^{-1}\left\langle \nabla_{\mathcal{A}^{\prime}} \mathcal{B}^{\prime
},\mathcal{R} \right\rangle ) \mid_{\Sigma}\;= (\tau^{-1}\left\langle
\nabla_{\mathcal{A}} \mathcal{B}, \mathcal{R} \right\rangle ) \mid_{\Sigma}\;
+ A(\tau) (\tau^{-1}\left\langle \bar{\mathcal{B}},\mathcal{R} \right\rangle
)\mid_{\Sigma}\;+ II(A,\bar{B}) + II(\bar{A},B) \;\;,\;\;
\]

\noindent and therefore%

\begin{equation}
\left.
\begin{array}
[c]{l}%
\nabla_{\mathcal{A}^{\prime}} \mathcal{B}^{\prime}\mid_{\Sigma}\;=
\nabla_{\mathcal{A}} \mathcal{B} \mid_{\Sigma}\; + \Sigma^{m-1}_{i=1}
\square_{A} (\mathcal{B}^{\prime}- \mathcal{B})(E_{i})\; E_{i} +\\
\\
\;\;\;\;\;\;\;\;\;\;\;\;\;\;\;\; + \{ A(\tau) (\tau^{-1}\left\langle
\bar{\mathcal{B}},\mathcal{R} \right\rangle )\mid_{\Sigma}\;+ II(A,\bar{B}) +
II(\bar{A},B) \}\;R \;\;.
\end{array}
\right.  \;\; \label{covar.deriv.restr.depends}%
\end{equation}

\bigskip

We have the following proposition, whose proof is straightforward using
(\ref{HII}) and (\ref{S}):

\smallskip

\begin{proposition}
\label{prop:covar.deriv.} Let $(M,g)$ be a transverse Riemann-Lorentz space
with tangent radical, singular hypersurface $\Sigma$ and canonical screen
distribution $S$. Let $\nabla$ be the Levi-Civita connection on $M-\Sigma$.

Let be $\mathcal{A},\mathcal{B}\in\frak{X}(M)$.

(a) It holds: $\;\nabla_{\mathcal{A}^{o}}\mathcal{B}^{o}\cong
0\;\;\Leftrightarrow\;\;II(A,B)=0$.

(b) The following two assertions are equivalent: (1) $\nabla_{\mathcal{A}^{o}%
}\mathcal{B}^{o}\cong0,$ whenever $\mathcal{A},\mathcal{B}$ are tangent to
$\Sigma$; and (2) $\Sigma\;$ is $II$-flat.

Let be $A,B\in\frak{X}_{\Sigma}\;$.

(c) $\;\nabla_{A}B$ is not well-defined, whenever one of the vectorfields
$A,B$ is either $N$ or $R$.

(d) The following two assertions are equivalent: (1) $\;\nabla_{A}B$ is
well-defined, whenever $\;A,B\in\Gamma(S)$; and (2) $\Sigma\,$ is $II$-flat
$\;\;\;\;\;\endproof$
\end{proposition}

Despite of Proposition \ref{prop:covar.deriv.}, some particular arrangements
with covariant derivatives are always "extension independent". We mention here
three cases:

\smallskip

(a) Let $\mathcal{A}$ be tangent to $\Sigma$ (i.e. $\mathcal{A} = \mathcal{X}%
$, with $X \in\frak{X}(\Sigma)$). Writing $\mathcal{X}^{\prime}=\mathcal{X} +
\tau\bar{\mathcal{X}},\;\mathcal{B}^{\prime}=\mathcal{B} + \tau\bar
{\mathcal{B}}$ (for some $\bar{\mathcal{X}},\bar{\mathcal{B}} \in\frak{X}%
(M)$), we obtain from (\ref{covar.deriv.restr.depends}):%

\begin{equation}
\nabla_{\mathcal{X}^{\prime}} \mathcal{B}^{\prime}\mid_{\Sigma}\;=
\nabla_{\mathcal{X}} \mathcal{B} \mid_{\Sigma}\; + \{II(X,\bar{B}) +
II(\bar{X},B) \}\;R\;\;,\;\;\;\; \label{covar.deriv.restr.depends,bis}%
\end{equation}

\noindent from which it follows%

\[
\left\{
\begin{array}
[c]{l}%
\left\langle \nabla_{\mathcal{X}^{\prime}}\mathcal{B}^{\prime}\mid_{\Sigma
}\;,C\right\rangle =\left\langle \nabla_{\mathcal{X}}\mathcal{B}\mid_{\Sigma
}\;,C\right\rangle \;\;,\;\;\forall C\in\frak{X}_{\Sigma}\\
II(\nabla_{\mathcal{X}^{\prime}}\mathcal{B}^{\prime}\mid_{\Sigma
}\;,Z)=II(\nabla_{\mathcal{X}}\mathcal{B}\mid_{\Sigma}\;,Z)\;\;,\;\;\forall
Z\in\frak{X}(\Sigma)\;\;.
\end{array}
\right.
\]

Last equation shows that, \emph{if $\;\Sigma\;$ is $II$-flat}, a natural map%

\begin{equation}
III:\frak{X}(\Sigma) \times\frak{X}(\Sigma) \times\frak{X}(\Sigma) \rightarrow
C^{\infty}(\Sigma)\;,\;(X,Y,Z) \mapsto II(\nabla_{\mathcal{X}} \mathcal{Y}
\mid_{\Sigma}\;,Z)\;\; \label{III}%
\end{equation}

\noindent arises, which turns out to be (straightforward computation)
\emph{$C^{\infty}(\Sigma)$-trilinear} and \emph{symmetric in its first two}
entries. Moreover, it is very easy to see that it holds:%

\[
III(\cdot, \cdot, R)\; = \mathcal{H} \;\;\in{\frak{T}^{o}_{2}}_{Sim} (\Sigma)
\;\;.
\]

\smallskip

Using (\ref{III}) and (\ref{HII}) we obtain (for all $\,X,Y \in\frak{X}%
(\Sigma)$ and $\,V \in\Gamma(S)$): $III(X,Y,V) = 0,\; III(V,R,R) = 0 $ and
$\;III(R,R,R) = -1.$

We say that $\Sigma$ is \emph{$III$-flat} if (it is $II$-flat and) it holds:%

\begin{equation}
III(V,W,R)=0 \;,\;\mathrm{for\;all}\;V,W \in\Gamma(S)\;\;; \label{IIIflat}%
\end{equation}

\noindent or, because of (\ref{HII}), if $\;\nabla_{\mathcal{V}}
\mathcal{W}\;$ is \emph{tangent to $\Sigma$}, for all $\mathcal{V},\mathcal{W}
\in\frak{X}(M)$ tangent to $S$. It follows from (\ref{IIflat}) that:
$\Sigma\;$ is $III$-flat if and only if it is $II$-flat and $\mathcal{H}$-flat.

\begin{remark}
If $\Sigma$ is $\,II$-flat but the radical is transverse (remember Remark
\ref{rem:IIflat}), a well-defined map $\;III:\frak{X}(\Sigma)\times
\frak{X}(\Sigma)\times\frak{X}_{\Sigma}\rightarrow C^{\infty}(\Sigma
)\;,\;(X,Y,C)\mapsto II(\nabla_{\mathcal{X}}\mathcal{Y}\mid_{\Sigma}\;,C)\,$
arises (see \cite{kossow2}; observe that the domain of $III$ in that case is
''larger'' than in our case with tangent radical), which is $C^{\infty}%
(\Sigma)$-trilinear and symmetric in the first two entries, and whose
restriction to $\,\frak{X}(\Sigma)\times\frak{X}(\Sigma)\times\frak{X}%
(\Sigma)\,$ vanishes. Then $\Sigma$ was defined \cite{kossow2} to be
$III$-flat if it holds:%

\[
III(X,Y,Rad)=0\;,\;for\;all\;\;X,Y\in\frak{X}(\Sigma)\;\;;
\]

\noindent although natural, this condition is in some sense stronger than
(\ref{IIIflat}) $\;\;\;\;\;\endproof$
\end{remark}

\smallskip

(b) Let $\mathcal{A}$ and $\mathcal{B}$ be tangent to $\Sigma$ (i.e.
$\mathcal{A}=\mathcal{X},\mathcal{B}=\mathcal{Y}$, with $X,Y \in
\frak{X}(\Sigma)$). Writing $\mathcal{X}^{\prime}=\mathcal{X} + \tau
\bar{\mathcal{X}},\;\mathcal{Y}^{\prime}=\mathcal{Y} + \tau\bar{\mathcal{Y}}$
(for some $\bar{\mathcal{X}},\bar{\mathcal{Y}} \in\frak{X}(M)$), we obtain
from the symmetry of $\nabla$ and (\ref{ro}):%

\begin{equation}
(\tau^{-1} \left\langle \mathcal{R},\nabla_{\mathcal{X}^{\prime}}
\mathcal{Y}^{\prime}\right\rangle _{Ant}) \mid_{\Sigma}\;= \frac{1}{2}
(\tau^{-1} \left\langle \mathcal{R},[\mathcal{X}^{\prime},\mathcal{Y}^{\prime
}] \right\rangle ) \mid_{\Sigma}\;= \frac{1}{2} \rho([X,Y])\;\;,\;\;
\label{nabla antysim. restr.}%
\end{equation}

\noindent where $\left\langle \mathcal{R},\nabla_{\mathcal{X}^{\prime}}
\mathcal{Y}^{\prime}\right\rangle _{Ant}$ means the antisymmetric part of
$\left\langle \mathcal{R},\nabla_{\mathcal{X}^{\prime}} \mathcal{Y}^{\prime
}\right\rangle $ under the permutation of $\mathcal{X}^{\prime}$ and
$\mathcal{Y}^{\prime}$.

\smallskip

(c) Finally, let be $X \in\frak{X}(\Sigma)$ and let $\mathcal{X},\mathcal{R}$
be the \emph{canonical} extensions of $X,R$. Let $\mathcal{X}^{\prime}$ be
\emph{any} extension of $X$. Writing $\mathcal{X}^{\prime}=\mathcal{X} +
\tau\bar{\mathcal{X}}$ (for some $\bar{\mathcal{X}} \in\frak{X}(M)$), we obtain:%

\begin{equation}
\left\{
\begin{array}
[c]{l}%
\nu(\nabla_{\mathcal{X}^{\prime}}\mathcal{R}\mid_{\Sigma}\;)=\square
_{X}R(N)=\mathcal{H}(X,R)\\
\\
\rho(\nabla_{\mathcal{X}^{\prime}}\mathcal{R}\mid_{\Sigma}\;)=(\tau
^{-1}\left\langle \mathcal{R},\nabla_{\mathcal{X}^{\prime}}\mathcal{R}%
\right\rangle )\mid_{\Sigma}\;=\frac{1}{2}(\tau^{-1}\mathcal{X}^{\prime}%
(\tau))\mid_{\Sigma}\;=\frac{1}{2}(\bar{\mathcal{X}}(\tau))\mid_{\Sigma}%
\;=\nu(\bar{X})\\
\\
\left\langle \nabla_{\mathcal{X}^{\prime}}\mathcal{R}\mid_{\Sigma
}\;,W\right\rangle =\square_{X}R(W)=-II(X,W)\;\;\;,\;\;\mathrm{for\;all}%
\;\;W\in\Gamma(S)
\end{array}
\right.  \label{ro,nu-nabla}%
\end{equation}

\noindent(in the middle line, we have used (\ref{ro}), (\ref{Xtau}) and
(\ref{Ntau}) in the first, third and last equality, respectively); in
particular: $\;\nabla_{\mathcal{R}}\mathcal{R}\mid_{\Sigma}\;=-N\;.$

\subsection{Covariant Curvatures}

Given $\mathcal{A},\mathcal{B},\mathcal{C},\mathcal{D} \in\frak{X}(M)$, we
compute Levi-Civita covariant curvature on $M -\Sigma$:%

\begin{equation}
\left.
\begin{array}
[c]{l}%
\left\langle R(\mathcal{A}^{o},\mathcal{B}^{o})\mathcal{C}^{o},\mathcal{D}^{o}
\right\rangle := \square_{\mathcal{A}^{o}} (\nabla_{\mathcal{B}^{o}}
\mathcal{C}^{o})(\mathcal{D}^{o}) - \square_{\mathcal{B}^{o}} (\nabla
_{\mathcal{A}^{o}} \mathcal{C}^{o})(\mathcal{D}^{o}) - \square_{[\mathcal{A}%
^{o},\mathcal{B}^{o}]} \mathcal{C}^{o}(\mathcal{D}^{o}) =\\
\\
\;\; = \sum^{m-1}_{i=1} \{ \mathcal{A}^{o} (\square_{\mathcal{B}^{o}}
\mathcal{C}^{o} (\mathcal{E}^{o}_{i})) - \mathcal{B}^{o} (\square
_{\mathcal{A}^{o}} \mathcal{C}^{o} (\mathcal{E}^{o}_{i})) \}\;\left\langle
\mathcal{E}^{o}_{i},\mathcal{D}^{o} \right\rangle +\\
\\
\;\;\;\;\; + \sum^{m-1}_{i=1} \{ \square_{\mathcal{B}^{o}} \mathcal{C}^{o}
(\mathcal{E}^{o}_{i}) \square_{\mathcal{A}^{o}} \mathcal{E}^{o}_{i}
(\mathcal{D}^{o}) - \square_{\mathcal{A}^{o}} \mathcal{C}^{o} (\mathcal{E}%
^{o}_{i}) \square_{\mathcal{B}^{o}} \mathcal{E}^{o}_{i} (\mathcal{D}^{o}) \}
+\\
\\
\;\;\;\;\; + \mathcal{A}^{o}(\square_{\mathcal{B}^{o}} \mathcal{C}^{o}
(\mathcal{R}^{o}) ({\tau^{o}}^{-1} \left\langle \mathcal{R},\mathcal{D}
\right\rangle ) - \mathcal{B}^{o}(\square_{\mathcal{A}^{o}} \mathcal{C}^{o}
(\mathcal{R}^{o}) ({\tau^{o}}^{-1} \left\langle \mathcal{R},\mathcal{D}
\right\rangle ) - \square_{[\mathcal{A}^{o},\mathcal{B}^{o}]} \mathcal{C}%
^{o}(\mathcal{D}^{o}) +\\
\\
\;\;\;\;\; + {\tau^{o}}^{-1} \; \{\square_{\mathcal{A}^{o}} \mathcal{C}^{o}
(\mathcal{R}^{o}) \square_{\mathcal{B}^{o}} \mathcal{D}^{o} (\mathcal{R}^{o})
- \square_{\mathcal{B}^{o}} \mathcal{C}^{o} (\mathcal{R}^{o}) \square
_{\mathcal{A}^{o}} \mathcal{D}^{o} (\mathcal{R}^{o}) \} \equiv\\
\\
\;\;\equiv\;{\tau^{o}}^{-1} \Upsilon(\mathcal{A}^{o},\mathcal{B}%
^{o},\mathcal{C}^{o},\mathcal{D}^{o}) \;\;,
\end{array}
\right.  \label{upsilon}%
\end{equation}

\noindent where the (everywhere regular) tensorfield $\;\Upsilon\in
\frak{T}^{0}_{4} (M)\;$ has the same symmetry properties as the covariant
curvature and satisfies (as usual, we denote $A \equiv\mathcal{A}\mid_{\Sigma
},... \in\frak{X}_{\Sigma}$):%

\begin{equation}
\left.
\begin{array}
[c]{l}%
\Upsilon(A,B,C,D)= \det\left(
\begin{array}
[c]{c}%
II(A,C) \;\; II(A,D)\\
II(B,C) \;\; II(B,D)
\end{array}
\right)  =\\
\\
\;\;\;\;\;\;\;\;\;\; \overset{\mathrm{If}\;\Sigma\;\mathrm{is}\;
II\mathrm{-flat}}{=} \;\;\det\left(
\begin{array}
[c]{c}%
\nu(A) \, \nu(B)\\
\rho(A) \, \rho(B)
\end{array}
\right)  \cdot\det\left(
\begin{array}
[c]{c}%
\nu(C) \, \nu(D)\\
\rho(C) \, \rho(D)
\end{array}
\right)  \;.
\end{array}
\right.  \label{upsilon restr.}%
\end{equation}

\smallskip

We obtain:

\smallskip

\begin{theorem}
\label{theor:covar.curvat.extends} Let $(M,g)$ be a transverse Riemann-Lorentz
space with tangent radical and singular hypersurface $\Sigma$. Let be
$\mathcal{A},\mathcal{B},\mathcal{C},\mathcal{D}\in\frak{X}(M)$.

(a) It holds: $\;\left\langle R(\mathcal{A}^{o},\mathcal{B}^{o})\mathcal{C}%
^{o},\mathcal{D}^{o}\right\rangle \cong0\;\;\Leftrightarrow\;\;\Upsilon(A,B,C,D)=0\;.$

(b) If we consider the assertions:

(1) $\;\Upsilon(A,B,C,D)=0\;\Leftarrow\;\det\left(
\begin{array}
[c]{c}%
\nu(A)\,\nu(B)\\
\rho(A)\,\rho(B)
\end{array}
\right)  \cdot\det\left(
\begin{array}
[c]{c}%
\nu(C)\,\nu(D)\\
\rho(C)\,\rho(D)
\end{array}
\right)  =0$

(2) $\;\Sigma\;\;is\;\;II-flat$

(3) $\;\Upsilon(A,B,C,D)=0\;\Rightarrow\;\det\left(
\begin{array}
[c]{c}%
\nu(A)\,\nu(B)\\
\rho(A)\,\rho(B)
\end{array}
\right)  \cdot\det\left(
\begin{array}
[c]{c}%
\nu(C)\,\nu(D)\\
\rho(C)\,\rho(D)
\end{array}
\right)  =0\;,$

\noindent then it holds: (1)$\,\Leftrightarrow\,$(2)$\,\Rightarrow\,$(3).
\end{theorem}

\begin{proof}
(a) follows from (\ref{HII}). (b) We compute $\Upsilon(A,B,C,D)$ from
(\ref{upsilon restr.}). To prove $\,(1)\Rightarrow\,(2),$ we first obtain:
$\;0=\Upsilon(V,R,W,N)=II(V,W),$ for all $V,W\in\Gamma(S),$ and the result
follows from (\ref{HII}) $\;\;\;\;\;$
\end{proof}

\begin{remark}
\label{rem:covar.curvat.extends} Theorem \ref{theor:covar.curvat.extends} also
says: $\;\Sigma$ is $II$-flat if and only if $\;\left\langle R(\mathcal{A}%
^{o},\mathcal{B}^{o})\mathcal{C}^{o},\mathcal{D}^{o}\right\rangle \cong0\,,$
whenever at least three of the vectorfields $\mathcal{A},\mathcal{B}%
,\mathcal{C},\mathcal{D}\in\frak{X}(M)$ are tangent to $\Sigma$ (as for the
sufficient condition: part (a) of the Theorem leads to $\;0=\Upsilon
(V,R,W,N)=II(V,W),$ for all $V,W\in\Gamma(S),$ and the result follows from
(\ref{HII})). Concerning other cases: $\;\left\langle R(\mathcal{N}%
^{o},\mathcal{V}^{o})\mathcal{N}^{o},\mathcal{W}^{o}\right\rangle \cong0\,$
(always) and $\;\left\langle R(\mathcal{N}^{o},\mathcal{V}^{o})\mathcal{N}%
^{o},\mathcal{R}^{o}\right\rangle \cong0\,$ (always), for all $\mathcal{V}%
,\mathcal{W}\in\frak{X}(M)$ tangent to $S$. However, even in the $II$-flat
case, $\;\left\langle R(\mathcal{N}^{o},\mathcal{R}^{o})\mathcal{N}%
^{o},\mathcal{R}^{o}\right\rangle \ncong0\,$.

These results contrast with the corresponding ones in the case of transverse
radical, namely (\cite{kossow2}, Theorem 3a): $\;\Sigma$ is $II$-flat if and
only if $\;\left\langle R(\mathcal{A}^{o},\mathcal{B}^{o})\mathcal{C}%
^{o},\mathcal{D}^{o}\right\rangle \cong0\,,$ for all vectorfields
$\;\mathcal{A},\mathcal{B},\mathcal{C},\mathcal{D}\in\frak{X}(M)$. Thus
$II$-flatness leads to no divergences in the case with transverse radical
$\;\;\;\;\;\endproof$
\end{remark}

Once $\left\langle R(\mathcal{A}^{o},\mathcal{B}^{o})\mathcal{C}%
^{o},\mathcal{D}^{o} \right\rangle \cong0 $, we may define (by continuity, the
only possible definition):%

\[
\left\langle R(\mathcal{A},\mathcal{B})\mathcal{C},\mathcal{D}\right\rangle
:=\tau^{-1}\Upsilon(\mathcal{A},\mathcal{B},\mathcal{C},\mathcal{D})\;\;,\;\;
\]

\noindent which satisfies (for all $\bar{\mathcal{A}} \in\frak{X}(M)$):%

\[
\left\langle R(\mathcal{A} + \tau\bar{\mathcal{A}}, \mathcal{B})\mathcal{C}%
,\mathcal{D} \right\rangle = \left\langle R(\mathcal{A}, \mathcal{B}%
)\mathcal{C},\mathcal{D} \right\rangle + \Upsilon(\bar{\mathcal{A}},
\mathcal{B},\mathcal{C},\mathcal{D}) \;.
\]

\begin{remark}
One should be careful in writing: $\left\langle R(\mathcal{A},\mathcal{B}%
)\mathcal{C},\mathcal{D}\right\rangle =\square_{\mathcal{A}}(\nabla
_{\mathcal{B}}\mathcal{C})(\mathcal{D})-\square_{\mathcal{B}}(\nabla
_{\mathcal{A}}\mathcal{C})(\mathcal{D})-\square_{\lbrack\mathcal{A}%
,\mathcal{B}]}\mathcal{C}(\mathcal{D}).$ The reason is that%

\[
\left\{
\begin{array}
[c]{l}%
\nabla_{\mathcal{B}^{o}}\mathcal{C}^{o}\cong{\tau^{o}}^{-1}\;\square
_{\mathcal{B}^{o}}\mathcal{C}^{o}(\mathcal{R}^{o})\;\mathcal{R}^{o}%
\;\;\;,\;\;whereas\\
\square_{\mathcal{A}^{o}}(\nabla_{\mathcal{B}^{o}}\mathcal{C}^{o}%
)(\mathcal{D}^{o})\cong-{\tau^{o}}^{-1}\;\square_{\mathcal{B}^{o}}%
\mathcal{C}^{o}(\mathcal{R}^{o})\square_{\mathcal{A}^{o}}\mathcal{D}%
^{o}(\mathcal{R}^{o})\;\;\;;
\end{array}
\right.
\]

\noindent thus it may happen that $\left\langle R(\mathcal{A},\mathcal{B}%
)\mathcal{C},\mathcal{D}\right\rangle $ is well-defined, but $\;\nabla
_{\mathcal{B}}\mathcal{D}$ does not exist (e.g. $\left\langle R(\mathcal{N}%
^{o},\mathcal{V}^{o})\mathcal{N}^{o},\mathcal{V}^{o}\right\rangle \cong0$, but
it may happen that $\;\nabla_{\mathcal{V}^{o}}\mathcal{V}^{o}\ncong
0)\;\;\endproof$
\end{remark}

As it happens with covariant derivatives, \emph{$\left\langle
R(A,B)C,D\right\rangle $ may sometimes be well-defined}. Indeed, starting with
$\mathcal{A}^{\prime}=\mathcal{A}+\tau\bar{\mathcal{A}},\;\mathcal{B}^{\prime
}=\mathcal{B}+\tau\bar{\mathcal{B}}\;\mathcal{C}^{\prime}=\mathcal{C}+\tau
\bar{\mathcal{C}}\;\mathcal{D}^{\prime}=\mathcal{D}+\tau\bar{\mathcal{D}}$
(for some $\bar{\mathcal{A}},\bar{\mathcal{B}},\bar{\mathcal{C}}%
,\bar{\mathcal{D}}\in\frak{X}(M)$), one sees that it holds (we denote $\bar
{A}\equiv\bar{\mathcal{A}}\mid_{\Sigma}\;,...$):%

\[
\left.
\begin{array}
[c]{l}%
\left\langle R(\mathcal{A}^{\prime},\mathcal{B}^{\prime})\mathcal{C}^{\prime
},\mathcal{D}^{\prime}\right\rangle \mid_{\Sigma}\;- \left\langle
R(\mathcal{A},\mathcal{B})\mathcal{C},\mathcal{D} \right\rangle \mid_{\Sigma
}\; =\\
\;\;\;\;\;\;\;\;\;\; = \Upsilon(\bar{A},B,C,D) + \Upsilon(A,\bar{B},C,D) +
\Upsilon(A,B,\bar{C},D) + \Upsilon(A,B,C,\bar{D}) \;\;.
\end{array}
\right.
\]

Last equation and (\ref{upsilon restr.}) lead to

\smallskip

\begin{theorem}
\label{theor:covar.curvat.restricts} Let $(M,g)$ be a transverse
Riemann-Lorentz space with tangent radical and singular hypersurface $\Sigma$.
Let be $A,B,C,D\in\frak{X}_{\Sigma}$. Then it holds:%

\[
\left.
\begin{array}
[c]{c}%
\left.
\begin{array}
[c]{c}%
\left\langle R(A,B)C,D\right\rangle \\
is\;well-defined
\end{array}
\right.  \;\;\Leftrightarrow\;\left.
\begin{array}
[c]{c}%
\Upsilon(\cdot,B,C,D)=\Upsilon(A,\cdot,C,D)=\;\;\;\;\;\\
\;\;\;\;\;=\Upsilon(A,B,\cdot,D)=\Upsilon(A,B,C,\cdot)=0
\end{array}
\right.  \;\Leftrightarrow\\
\\
\overset{If\;\;\Sigma\;\;is\;II-flat}{\Leftrightarrow}\;\;\left\{
\begin{array}
[c]{l}%
either\;\;\det\left(
\begin{array}
[c]{c}%
\nu(A)\;\nu(B)\\
\rho(A)\;\rho(B)
\end{array}
\right)  =0=\det\left(
\begin{array}
[c]{c}%
\nu(C)\;\nu(D)\\
\rho(C)\;\rho(D)
\end{array}
\right)  \\
or\;\;one\;of\;the\;above\;two\;matrices\;vanishes
\end{array}
\right.  \;\endproof
\end{array}
\right.
\]
\end{theorem}

Let be $p\in\Sigma\,$ and $\,a,b,c,d\in T_{p}M$. We say that
\emph{$\left\langle R(a,b)c,d\right\rangle $ is well-defined} if there exist
(local) extensions $\,A,B,C,D\in\frak{X}_{\Sigma}\;$ of $\,a,b,c,d$ such that
$\left\langle R(A,B)C,D\right\rangle $ is well-defined. This definition is
independent of the chosen extensions $A,B,C,D$, as the next Lemma shows:

\smallskip

\begin{lemma}
\label{covar.curvat.pointwise} Let $(M,g)$ be a transverse Riemann-Lorentz
space with tangent radical and singular hypersurface $\Sigma$.

(a) Let be $\,p\in\Sigma$ and $\,a,b,c,d\in T_{p}M$ such that $\left\langle
R(a,b)c,d\right\rangle $ is well-defined. Let $A,B,C,D\in\frak{X}_{\Sigma}\;$
and $A^{\prime},B^{\prime},C^{\prime},D^{\prime}\in\frak{X}_{\Sigma}\;$ be two
sets of (local) extensions of $\,a,b,c,d$ such that both $\left\langle
R(A,B)C,D\right\rangle $ and $\left\langle R(A^{\prime},B^{\prime})C^{\prime
},D^{\prime}\right\rangle $ are well-defined. Then it holds:%

\[
\left\langle R(A^{\prime},B^{\prime})C^{\prime},D^{\prime}\right\rangle
(p)=\left\langle R(A,B)C,D\right\rangle (p)\;\;.
\]

$(b)$ Let be $A,B,C,D\in\frak{X}_{\Sigma}\;$ such that $\;\left\langle
R(A(p),B(p))C(p),D(p)\right\rangle $ is well-defined, for all $p\in\Sigma$.
Then $\left\langle R(A,B)C,D\right\rangle $ is well-defined and it holds (for
all $p\in\Sigma$):%

\[
\left\langle R(A,B)C,D\right\rangle (p)=\left\langle
R(A(p),B(p))C(p),D(p)\right\rangle \;\;.
\]
\end{lemma}

\begin{proof}
(a) Let $\mathcal{A},\mathcal{B},\mathcal{C},\mathcal{D}\in\frak{X}(M)$ and
$\mathcal{A}^{\prime},\mathcal{B}^{\prime},\mathcal{C}^{\prime},\mathcal{D}%
^{\prime}\in\frak{X}(M)$ be the \emph{canonical} extensions of $A,B,C,D$ and
$A^{\prime},B^{\prime},C^{\prime},D^{\prime}$, respectively. Let us consider
the function $F:=\left\langle R(\mathcal{A}^{\prime},\mathcal{B}^{\prime
})\mathcal{C}^{\prime},\mathcal{D}^{\prime}\right\rangle -\left\langle
R(\mathcal{A},\mathcal{B})\mathcal{C},\mathcal{D}\right\rangle \,\in
C^{\infty}(M),$ which satisfies (on $M-\Sigma$)%

\[
\left.
\begin{array}
[c]{l}%
F^{o}=\left\langle R({\mathcal{A}^{\prime}}^{o}-\mathcal{A}^{o},{\mathcal{B}%
^{\prime}}^{o}){\mathcal{C}^{\prime}}^{o},{\mathcal{D}^{\prime}}%
^{o}\right\rangle +\left\langle R(\mathcal{A}^{o},{\mathcal{B}^{\prime}}%
^{o}-\mathcal{B}^{o}){\mathcal{C}^{\prime}}^{o},{\mathcal{D}^{\prime}}%
^{o}\right\rangle +\\
\;\;\;\;\;\;\;\;\;\;\;\;\;\;\;\;\;\;\;\;+\left\langle (\mathcal{A}%
^{o},\mathcal{B}^{o})({\mathcal{C}^{\prime}}^{o}-\mathcal{C}^{o}%
),{\mathcal{D}^{\prime}}^{o}\right\rangle +\left\langle R(\mathcal{A}%
^{o},\mathcal{B}^{o})\mathcal{C}^{o},{\mathcal{D}^{\prime}}^{o}-\mathcal{D}%
^{o}\right\rangle \;.
\end{array}
\right.
\]

Let $\gamma$ be the integral curve of the canonical extension $\mathcal{N}$ of
$N$ with $\gamma(0)=p$. Because $\,\mathcal{A}^{\prime}-\mathcal{A}\,$ is the
canonical extension of $\,A^{\prime}-A\,$ and $(A^{\prime}-A)(p)=0$, it
follows that $(\mathcal{A}^{\prime}-\mathcal{A})\circ\gamma=0$. Analogously,
$\,(\mathcal{B}^{\prime}-\mathcal{B})\circ\gamma=(\mathcal{C}^{\prime
}-\mathcal{C})\circ\gamma=(\mathcal{D}^{\prime}-\mathcal{D})\circ\gamma=0\,$.
Therefore $\;F^{o}\circ\gamma=0\;$, and we obtain%

\[
\left.
\begin{array}
[c]{l}%
\left\langle R(A^{\prime},B^{\prime})C^{\prime},D^{\prime}\right\rangle
(p)-\left\langle R(A,B)C,D\right\rangle (p)=\;\;(\mathrm{by\;hypothesis})\\
\\
\;\;\;\;\;=\left\langle R(\mathcal{A}^{\prime},\mathcal{B}^{\prime
})\mathcal{C}^{\prime},\mathcal{D}^{\prime}\right\rangle (p)-\left\langle
R(\mathcal{A},\mathcal{B})\mathcal{C},\mathcal{D}\right\rangle (p)=\lim
_{t\rightarrow0}(F^{o}\circ\gamma)(t)=0\;.
\end{array}
\right.
\]

(b) follows immediately from Theorem \ref{theor:covar.curvat.restricts},
because the condition there is tensorial $\;\;\;\;\;$
\end{proof}

\bigskip

Theorem \ref{theor:covar.curvat.restricts} also gives the algebraic (necessary
and sufficient) condition on $\,a,b,c,d\;$ in order that $\;\left\langle
R(a,b)c,d\right\rangle $ is well-defined, namely: $\;\Upsilon(\cdot
,b,c,d)=\Upsilon(a,\cdot,c,d)=\Upsilon(a,b,\cdot,d)=\Upsilon(a,b,c,\cdot)=0$.

\subsection{Sectional Curvatures}

We start with two $C^{\infty}(M)$-linearly independent vectorfields
$\mathcal{A},\mathcal{B} \in\frak{X}(M)$ with $\;rank (g_{\mathcal{A}^{o}
\wedge\mathcal{B}^{o}}) = 2$ and compute the $\nabla$-sectional curvature on
$M -\Sigma$%

\[
K_{\mathcal{A}^{o} \wedge\mathcal{B}^{o}} := \frac{\left\langle R(\mathcal{A}%
^{o},\mathcal{B}^{o})\mathcal{A}^{o}, \mathcal{B}^{o} \right\rangle }
{\det(g_{\mathcal{A}^{o} \wedge\mathcal{B}^{o}})} \;\;.
\]

\bigskip

(a) \emph{Suppose first that $\;rank (g_{A \wedge B}) = \;(constant)\; 2$} (an
open condition). In that case, we have:

\smallskip

\begin{proposition}
\label{prop:sec.curvat.extends1} Let $(M,g)$ be a transverse Riemann-Lorentz
space with tangent radical and singular hypersurface $\Sigma$. Let
$\mathcal{A},\mathcal{B}\in\frak{X}(M)$ be such that $rank(g_{\mathcal{A}%
\wedge\mathcal{B}})=2$.

(a) It holds: $\;K_{\mathcal{A}^{o}\wedge\mathcal{B}^{o}}\cong
0\;\;\;\Leftrightarrow\;\;\;\;\;\;\Upsilon(A,B,A,B)=0\;\;.$

(b) If we consider the assertions:

(1) $\;\Upsilon(A,B,A,B)=0\;\Leftarrow\;\det\left(
\begin{array}
[c]{c}%
\nu(A)\,\nu(B)\\
\rho(A)\,\rho(B)
\end{array}
\right)  =0$

(2) $\;\Sigma\;\;is\;II-flat$

(3) $\;\Upsilon(A,B,A,B)=0\;\Rightarrow\;\det\left(
\begin{array}
[c]{c}%
\nu(A)\,\nu(B)\\
\rho(A)\,\rho(B)
\end{array}
\right)  =0\;,$

\noindent then it holds: (1)$\,\Leftrightarrow\,$(2)$\,\Rightarrow\,$(3).
\end{proposition}

\begin{proof}
(a) follows from Theorem \ref{theor:covar.curvat.extends}a. (b) follows from
Theorem \ref{theor:covar.curvat.extends}b. To prove $\,(1)\Rightarrow\,(2),$
choose $A\in\Gamma(S)$ and $B=\nu(B)N+\rho(B)R,$ with both $\nu(B),\rho(B)\in
C^{\infty}(\Sigma)$ nowhere zero $\;\;\;\;\;$
\end{proof}

\begin{remark}
\label{rem:sec.curvat.extends1} It follows from Proposition
\ref{prop:sec.curvat.extends1} that: $\;K_{\mathcal{N}^{o}\wedge
\mathcal{V}^{o}}\cong0$ (always) and $\;K_{\mathcal{N}^{o}\wedge
(\mathcal{V}^{o}+\mathcal{R}^{o})}\ncong0$ (always), for all $\mathcal{V}%
\in\frak{X}(M)$ nowhere zero and tangent to $S$.

Moreover, if $\Sigma$ is $II$-flat, we obtain: $\;K_{\mathcal{A}^{o}%
\wedge\mathcal{B}^{o}}\cong0$, whenever both vectorfields $\mathcal{A}%
,\mathcal{B}\in\frak{X}(M)$ (with $rank(g_{\mathcal{A}\wedge\mathcal{B}})=2$)
are tangent to $\Sigma$, and $\;K_{(\mathcal{N}^{o}+\mathcal{R}^{o}%
)\wedge\mathcal{V}^{o}}\cong0$, for all $\mathcal{V}\in\frak{X}(M)$ nowhere
zero and tangent to $S$.

Things are again different when the radical is transverse; for then,
$II$-flatness is enough to guarantee (\cite{kossow2}, Theorem 3b) that
$\;K_{\mathcal{A}^{o}\wedge\mathcal{B}^{o}}\cong0$, for all $\mathcal{A}%
,\mathcal{B}\in\frak{X}(M)$ (with $rank(g_{\mathcal{A}\wedge\mathcal{B}})=2$)
$\;\;\;\;\;\endproof$
\end{remark}

Once $K_{\mathcal{A}^{o} \wedge\mathcal{B}^{o}} \cong0 $, we define:%

\[
K_{\mathcal{A}\wedge\mathcal{B}} := \frac{\left\langle R(\mathcal{A}%
,\mathcal{B})\mathcal{A}, \mathcal{B} \right\rangle } {\det(g_{\mathcal{A}
\wedge\mathcal{B}})} \;\;.
\]

Now we obtain:

\begin{proposition}
\label{prop:sec.curvat.restricts1} Let $(M,g)$ be a transverse Riemann-Lorentz
space with tangent radical and singular hypersurface $\Sigma$. Let
$A,B\in\frak{X}_{\Sigma}\;$ be such that $\;rank(g_{A\wedge B})=2$.

(a) It holds:%

\[
\left.
\begin{array}
[c]{c}%
K_{A\wedge B}\;\;is\;well-defined\;\;\;\Leftrightarrow\;\;\;\Upsilon
(\cdot,B,A,B)=\Upsilon(A,\cdot,A,B)=0\;\;\;\Leftrightarrow\\
\\
\overset{\mathrm{If}\;\Sigma\;\mathrm{is}\;II-\mathrm{flat}}{\Leftrightarrow
}\;\;\det\left(
\begin{array}
[c]{c}%
\nu(A)\;\nu(B)\\
\rho(A)\;\rho(B)
\end{array}
\right)  =0\;\;\;(\;\Leftrightarrow\;dim\;Span(A^{S},B^{S})\geq1\;)\;\;.
\end{array}
\right.
\]

(b) The following three assertions are equivalent:

\smallskip

(1) $\;K_{A\wedge B}\;\;is\;well-defined\;\;\;\Leftarrow\;\;\;dim\;Span(A^{S}%
,B^{S})\geq1$

(2) $\;\Sigma\;\;is\;II-flat$

(3) $\;K_{A\wedge B}\;\;is\;well-defined\;\;\;\Leftarrow\;\;\;dim\;Span(A^{S}%
,B^{S})=1$
\end{proposition}

\begin{proof}
(a) is a direct consequence of Theorem \ref{theor:covar.curvat.restricts}.

(b) We prove $\,(1)\Rightarrow(2)\,$: if $\Sigma$ is not $II$-flat, then
$\exists V\in\Gamma(S)$ such that $II(V,V)\neq0,$ thus (\ref{upsilon restr.})
leads to $\,\Upsilon(N,V,V,R)\neq0$ and Theorem
\ref{theor:covar.curvat.extends}a shows that $\left\langle R(\mathcal{N}%
^{o},\mathcal{V}^{o})\mathcal{V}^{o},\mathcal{R}^{o}\right\rangle \ncong0$
(for all extensions $\mathcal{N},\mathcal{V},\mathcal{R}$ of $N,V,R$). But
this cannot be true, because it holds:%

\[
\left.
\begin{array}
[c]{l}%
\left\langle R(\mathcal{N}^{o}+\mathcal{R}^{o},\mathcal{V}^{o})(\mathcal{N}%
^{o}+\mathcal{R}^{o}),\mathcal{V}^{o}\right\rangle =\left\langle
R(\mathcal{N}^{o},\mathcal{V}^{o})\mathcal{N}^{o},\mathcal{V}^{o}\right\rangle
+\\
\;\;\;\;\;\;\;\;\;\;\;\;\;\;\;\;\;\;\;\;\;\;\;\;\;\;\;\;\;\;+\left\langle
R(\mathcal{R}^{o},\mathcal{V}^{o})\mathcal{R}^{o},\mathcal{V}^{o}\right\rangle
+2\left\langle R(\mathcal{N}^{o},\mathcal{V}^{o})\mathcal{R}^{o}%
,\mathcal{V}^{o}\right\rangle \;,
\end{array}
\right.
\]

\noindent and the first three terms have good limits on $\Sigma$ because of
Theorem \ref{theor:covar.curvat.extends}a.

Now $\,(2)\Rightarrow(3)\,$ follows from Part (a). Finally, $\,(3)\Rightarrow
(1)\,$ is a consequence of Proposition \ref{prop:sec.curvat.extends1} $\;\;\;\;\;$
\end{proof}

\bigskip

Let be $\;p\in\Sigma\;$ and $\;a,b\in T_{p}M$ such that $\;rank(g_{a\wedge
b})=2$. We say that \emph{$K_{a\wedge b}$ is well-defined} if there exist
(local) extensions $A,B\in\frak{X}_{\Sigma}$ of $a,b$ such that $K_{A\wedge
B}$ is well-defined\footnote{Such extensions $A,B$ will always satisfy,
locally around $p$, $\;rank(g_{A\wedge B})=2$}. This definition is independent
of the chosen extensions $A,B$, as the above Lemma shows. Proposition
\ref{prop:sec.curvat.restricts1}a also gives the algebraic (necessary and
sufficient) condition on $\,a,b$ in order that $K_{a\wedge b}$ is
well-defined, namely: $\;\Upsilon(\cdot,b,a,b)=\Upsilon(a,\cdot,a,b)=0$.

\bigskip

(b) \emph{Suppose now that $\;rank (g_{A \wedge B}) <2$} $\;(\,\Leftrightarrow
\;R \in A\wedge B \;,\; \Leftrightarrow rank(g_{A\wedge B})=1\,)$. It follows:%

\begin{equation}
\left.
\begin{array}
[c]{l}%
\det(g_{\mathcal{A} \wedge\mathcal{B}}) = k \tau\;,\;k \in C^{\infty}(M),
\;\mathrm{with}\; k\mid_{\Sigma}\;\neq0 \;(\mathrm{everywhere})
\;;\;\Rightarrow\\
\\
\;\;\;\;\;\;\;\;\;\; \Rightarrow\;\;(\nu(A)N+A^{S}) \wedge(\nu(B)N+B^{S}) = 0
\end{array}
\right.  \label{rank1}%
\end{equation}

From (\ref{upsilon restr.}) we compute $\;\Upsilon(A,B,A,B) = - \det^{2}
\left(
\begin{array}
[c]{c}%
\nu(A) \; \nu(B)\\
\rho(A) \; \rho(B)
\end{array}
\right)  \;$ and, using (\ref{rank1}), we conclude:%

\begin{equation}
\Upsilon(A,B,A,B)=0 \;\;\Leftrightarrow\;\;\nu(A)=0=\nu(B) \;\;.
\label{rank1bis}%
\end{equation}

We finally obtain:

\begin{proposition}
\label{prop:sec.curvat.extends2} Let $(M,g)$ be a transverse Riemann-Lorentz
space with tangent radical and singular hypersurface $\Sigma$. Let
$\mathcal{A},\mathcal{B}\in\frak{X}(M)$ be such that $\;rank(g_{\mathcal{A}%
^{o}\wedge\mathcal{B}^{o}})=2,\;rank(g_{A\wedge B})=1$ and $\;R\in A\wedge B$.

(a) It holds: $K_{\mathcal{A}^{o}\wedge\mathcal{B}^{o}}\cong0\,\Leftrightarrow
\,\left\langle R(\mathcal{A},\mathcal{B})\mathcal{A},\mathcal{B}\right\rangle
\mid_{\Sigma}\;=0\,\Rightarrow\,\nu(A)=0=\nu(B)\,.$

(b) If we consider the assertions:

(1) $\left\langle R(A,B)A,B\right\rangle =0\Leftarrow\nu(A)=0=\nu
(B)\,(\Leftrightarrow dim\,Span(A^{S},B^{S})=1)$

(2) $\;\Sigma\;\;is\;III-flat$

(3) $\;\Sigma\;\;is\;II-flat$

(4) $\;\left\langle R(A,B)A,B\right\rangle \;\;is\;well-defined\;\;\Leftarrow
\;\;\nu(A)=0=\nu(B)\;\;,$

\noindent then it holds: (1)$\,\Leftrightarrow\,$(2)$\,\Rightarrow
\,$(3)$\,\Leftrightarrow\,$(4).
\end{proposition}

\begin{proof}
(a) follows from (\ref{rank1}), Theorem \ref{theor:covar.curvat.extends}a and
(\ref{rank1bis}).

(b) When $\nu(A)=0=\nu(B)$, it can be easily proved (as a consequence of
Theorem \ref{theor:covar.curvat.restricts}) that:%

\begin{equation}
\left\langle R(A,B)A,B\right\rangle \;\mathrm{is\;well-defined}%
\;\;\Leftrightarrow\;\;II(A^{S},A^{S})=0=II(B^{S},B^{S}%
)\;.\;\;\label{rank1tri}%
\end{equation}

Thus (3)$\,\Rightarrow\,$(4) becomes trivial. To prove (4)$\,\Rightarrow
\,$(3): if $\Sigma$ is not $II$-flat, then $\exists V\in\Gamma(S)$ such that
$II(V,V)\neq0$, and (\ref{rank1tri}) shows that $\left\langle
R(V,R)V,R\right\rangle $ cannot be well-defined.

We have (trivially): (1)$\,\Rightarrow\,$(4) and (2)$\,\Rightarrow\,$(3).

Because $\;R\in A\wedge B$, to prove (1)$\,\Leftrightarrow\,$(2) is equivalent
to prove that $\Sigma$ is $III$-flat if and only if $\;\left\langle
R(V,R)V,R\right\rangle $ is well-defined and $\,=0,$ for all $V\in\Gamma(S)$.
But it holds:%

\[
\left\langle R(\mathcal{V}^{o},\mathcal{R}^{o})\mathcal{V}^{o},\mathcal{R}%
^{o}\right\rangle =\;\square_{\mathcal{V}^{o}}(\nabla_{\mathcal{R}^{o}%
}\mathcal{V}^{o})(\mathcal{R}^{o})-\square_{\mathcal{R}^{o}}(\nabla
_{\mathcal{V}^{o}}\mathcal{V}^{o})(\mathcal{R}^{o})-\square_{\lbrack
\mathcal{V}^{o},\mathcal{R}^{o}]}\mathcal{V}^{o}(\mathcal{R}^{o})\;.
\]

As we have already seen, it follows from either $(1)$ or $(2)$ that $\Sigma$
is $II$-flat. On the other hand, $\;\nu(\nabla_{\mathcal{R}}\mathcal{V}%
\mid_{\Sigma})=\nu(\nabla_{\mathcal{V}}\mathcal{R}\mid_{\Sigma})=\square
_{V}R(N)=:\mathcal{H}(V,R)=0$ (thus $\nabla_{\mathcal{R}}\mathcal{V}$ is
\emph{tangent to $\Sigma$}). It thus follows that the first and third terms in
the righthand side are proportional to $\tau^{o}$. And we obtain:%

\[
\left.
\begin{array}
[c]{c}%
\underset{\forall\mathcal{V}\;\mathrm{tangent\;to}\;S}{\left\langle
R(\mathcal{V}^{o},\mathcal{R}^{o})\mathcal{V}^{o},\mathcal{R}^{o}\right\rangle
=k_{1}^{o}\tau^{o}}\;\;\Leftrightarrow\;\;\underset{\forall\mathcal{V}%
\;\mathrm{tangent\;to}\;S}{\square_{\mathcal{R}^{o}}(\nabla_{\mathcal{V}^{o}%
}\mathcal{V}^{o})(\mathcal{R}^{o})=k_{2}^{o}\tau^{o}}\;\Leftrightarrow\\
\\
\Leftrightarrow\;\;\underset{\forall\mathcal{V}\;\mathrm{tangent\;to}%
\;S}{II(\nabla_{\mathcal{V}}\mathcal{V}\mid_{\Sigma}\;,R)=0}%
\;\;\Leftrightarrow\;\;\underset{\forall V\in\Gamma(S)}{III(V,V,R)=0}%
\;\;\Leftrightarrow\\
\\
\Leftrightarrow\;\;\;\underset{\forall V,W\in\Gamma(S)}{III(V,W,R)=0}%
\;\;\;\Leftrightarrow\;\;\underset{\forall V\in\Gamma(S)\;\mathrm{and}%
\;\forall T\in\frak{X}(\Sigma)}{III(V,T,R)=0}\;\;\Leftrightarrow
\;\;\Sigma\;\mathrm{is}\;III\mathrm{-flat}\;;
\end{array}
\right.
\]

\noindent in the last five implications we have used Proposition
\ref{prop:covar.deriv.}a, (\ref{III}), the symmetry of $III$, the fact that
$III(V,R,R)=0$ and $II$-flatness, respectively $\;\;\;\;\;$
\end{proof}

\begin{remark}
It follows from Proposition \ref{prop:sec.curvat.extends2} that:
$\;K_{\mathcal{N}^{o}\wedge\mathcal{R}^{o}}\ncong0$ (always).

Moreover, if $\Sigma$ is $II$-flat but not $III$-flat, there exists
$\,V\in\Gamma(S)$ such that $\;\mathcal{H}(V,V)\neq0$. Let $\mathcal{V}%
\in\frak{X}(M)$ be \emph{any} extension of $V$. Using (\ref{upsilon}) and
$II$-flatness, we explicitly compute%

\[
\left\langle R(\mathcal{V}^{o},\mathcal{R}^{o})\mathcal{V}^{o},\mathcal{R}%
^{o}\right\rangle =k^{o}\tau^{o}-\square_{\mathcal{V}^{o}}\mathcal{V}%
^{o}(\mathcal{E}_{1}^{o})\square_{\mathcal{R}^{o}}\mathcal{E}_{1}%
^{o}(\mathcal{R}^{o})
\]

\noindent($k\in C^{\infty}(M)$). Because $\;II(R,N)=1,$ it follows:
$\;\left\langle R(\mathcal{V},\mathcal{R})\mathcal{V},\mathcal{R}\right\rangle
\mid_{\Sigma}\;=-\mathcal{H}(V,V)\neq0\;$ and, by Theorem
\ref{theor:covar.curvat.extends}a, $\;K_{\mathcal{V}^{o}\wedge\mathcal{R}^{o}%
}\ncong0$.

This result again contrasts with the corresponding one in the case of
transverse radical, where $III$-flatness guarantees (\cite{kossow2}, Theorem
3b) that $\;K_{\mathcal{A}^{o}\wedge\mathcal{B}^{o}}\cong0$, for all
$\mathcal{A},\mathcal{B}\in\frak{X}(M)$ (with $rank(g_{\mathcal{A}%
\wedge\mathcal{B}})=1$) $\;\;\;\;\;\endproof$
\end{remark}

Suppose that $K_{\mathcal{A}^{o}\wedge\mathcal{B}^{o}}\cong0$. Because, in
this case, both $\;\tau^{-1}\left\langle R(\mathcal{A},\mathcal{B}%
)\mathcal{A},\mathcal{B}\right\rangle \;$ (by Proposition
\ref{prop:sec.curvat.extends2}a) and $\;\tau^{-1}\det(g_{\mathcal{A}%
\wedge\mathcal{B}})\;$ (by (\ref{rank1})) must be well-defined functions and
the second one nowhere vanishes, we define:%

\[
K_{\mathcal{A} \wedge\mathcal{B}} := \frac{\tau^{-1}\left\langle
R(\mathcal{A},\mathcal{B})\mathcal{A},\mathcal{B} \right\rangle } {\tau
^{-1}\det(g_{\mathcal{A} \wedge\mathcal{B}})}\;\;.
\]

\bigskip

And finally we obtain:

\begin{proposition}
\label{prop:sec.curvat.restricts2} Let $(M,g)$ be a transverse Riemann-Lorentz
space with tangent radical and singular hypersurface $\Sigma$. Let
$A,B\in\frak{X}_{\Sigma}\;$ be such that $\;rank(g_{A\wedge B})=1$ and $\;R\in
A\wedge B$. Then: $K_{A\wedge B}$ is not well-defined.
\end{proposition}

\begin{proof}
In order to have $K_{\mathcal{A}^{o}\wedge\mathcal{B}^{o}}\cong0$, we can
write (Proposition \ref{prop:sec.curvat.extends2}a) $A=V\in\Gamma(S)$ and
$B=R$. Let us consider extensions $\mathcal{V},\mathcal{V}^{\prime
}=\mathcal{V}+\tau\,f\,\mathcal{N}$ of $V$ and $\mathcal{R}$ of $R$ (with
$\mathcal{N},\mathcal{V},\mathcal{R}$ \emph{canonical} and $f\in C^{\infty
}(M))$. Then it holds (by Theorem \ref{theor:covar.curvat.extends}a, all terms
in the next equality are well-defined):%

\[
\left.
\begin{array}
[c]{l}%
\left\langle R(\mathcal{V}^{\prime},\mathcal{R})\mathcal{V}^{\prime
},\mathcal{R}\right\rangle =\left\langle R(\mathcal{V},\mathcal{R}%
)\mathcal{V},\mathcal{R}\right\rangle +\\
\;\;\;\;\;\;\;\;\;\;\;\;\;\;\;\;\;\;\;\;+2\tau f\left\langle R(\mathcal{N}%
,\mathcal{R})\mathcal{V},\mathcal{R}\right\rangle +\tau f^{2}\Upsilon
(\mathcal{N},\mathcal{R},\mathcal{N},\mathcal{R})\;;\;\Rightarrow\\
\\
\Rightarrow\;\;(\tau^{-1}\left\langle R(\mathcal{V}^{\prime},\mathcal{R}%
)\mathcal{V}^{\prime},\mathcal{R}\right\rangle )\mid_{\Sigma}\;=(\tau
^{-1}\left\langle R(\mathcal{V},\mathcal{R})\mathcal{V},\mathcal{R}%
\right\rangle )\mid_{\Sigma}\;+\\
\;\;\;\;\;\;\;\;\;\;\;\;\;\;\;\;\;\;\;\;+2f\mid_{\Sigma}\;\left\langle
R(\mathcal{N},\mathcal{R})\mathcal{V},\mathcal{R}\right\rangle \mid_{\Sigma
}\;+f^{2}\mid_{\Sigma}\Upsilon(N,R,N,R)\;.
\end{array}
\right.
\]

But it follows from (\ref{upsilon restr.}) that $\;\Upsilon(N,R,N,R)=-1$;
therefore a suitable choice of the function $\,f\,$ leads, again by
Proposition \ref{prop:sec.curvat.extends2}a, to the result: $\;K_{\mathcal{V}%
^{\prime}\wedge\mathcal{R}}\mid_{\Sigma}\;\neq\;K_{\mathcal{V}\wedge
\mathcal{R}}\mid_{\Sigma}\;\;\;$
\end{proof}

\bigskip

Let be $\;p\in\Sigma\;$ and $\;a,b\in T_{p}M$ such that $\;rank(g_{a\wedge
b})<2$ $\;(\,\Leftrightarrow\;R(p)\in a\wedge b\;,\;\Leftrightarrow
rank(g_{a\wedge b})=1\,)$. We say that \emph{$K_{a\wedge b}$ is well-defined}
if there exist (local) extensions $A,B\in\frak{X}_{\Sigma}$ of $a,b$ with
$\;rank(g_{A\wedge B})<2$ and such that $K_{A\wedge B}$ is well-defined. By
Proposition \ref{prop:sec.curvat.restricts2}, \emph{$K_{a\wedge b}$ cannot be
well-defined}.

\subsection{Ricci Curvatures}

We start with two vectorfields $\mathcal{A},\mathcal{B} \in\frak{X}(M)$ and
compute Ricci-curvature on $M -\Sigma\;$:%

\begin{equation}
Ric(\mathcal{A}^{o},\mathcal{B}^{o}) = \sum^{m-1}_{i=1} \left\langle
R(\mathcal{A}^{o},\mathcal{E}^{o}_{i})\mathcal{B}^{o}, \mathcal{A}^{o}_{i}
\right\rangle + {\tau^{o}}^{-1} \;\left\langle R(\mathcal{A}^{o}%
,\mathcal{R}^{o})\mathcal{B}^{o}, \mathcal{R}^{o} \right\rangle \;\;.
\label{ricci}%
\end{equation}

Therefore: $Ric(\mathcal{A}^{o},\mathcal{B}^{o})\cong0$, whenever all $m$
terms on the right hand side have good limits on $\Sigma$. As for the last
term, this is equivalent (see \ref{notations}) to say that $\;\left\langle
R(\mathcal{A},\mathcal{R})\mathcal{B},\mathcal{R}\right\rangle \mid_{\Sigma
}\;=0$.

\smallskip

Now we have:

\begin{proposition}
\label{prop:ricci curvat.extends} Let $(M,g)$ be a transverse Riemann-Lorentz
space with tangent radical and singular hypersurface $\Sigma$. Given any
extensions $\mathcal{N}^{\prime}=\mathcal{N}+\tau\bar{\mathcal{N}}%
,\mathcal{V}^{\prime}=\mathcal{V}+\tau\bar{\mathcal{V}},\mathcal{R}^{\prime
}=\mathcal{R}+\tau\bar{\mathcal{R}}$ of $N,V\in\Gamma(S),R$ (where
$\mathcal{N},\mathcal{V},\mathcal{R}$ are the canonical extensions and
$\bar{\mathcal{N}},\bar{\mathcal{V}},\bar{\mathcal{R}}\in\frak{X}(M)$), it holds:

(a) $\;Ric({\mathcal{N}^{\prime}}^{o},{\mathcal{N}^{\prime}}^{o})\ncong0$

(b) $\;Ric({\mathcal{N}^{\prime}}^{o},{\mathcal{V}^{\prime}}^{o}%
)\cong0\;\;\Leftrightarrow\;\;\rho([V,R])+\nu(\bar{V})=0$

(c) $\;\Sigma\;\;is\;II-flat\;\;\Rightarrow\;\;\;[\;Ric({\mathcal{N}^{\prime}%
}^{o},{\mathcal{R}^{\prime}}^{o})\cong0\;\;\Leftrightarrow\;\;\nu(\bar{R})=0\;]$

(d) $\;Ric({\mathcal{R}^{\prime}}^{o},{\mathcal{V}^{\prime}}^{o}%
)\cong0\;\;,\;\;for\;all\;V\in\Gamma(S)$

(e) $\;Ric({\mathcal{R}^{\prime}}^{o},{\mathcal{R}^{\prime}}^{o})\ncong0$

(f) $\;\Sigma\;\mathrm{is}\;II-flat\;\;\Rightarrow\;\;[\;Ric({\mathcal{V}%
^{\prime}}^{o},{\mathcal{W}^{\prime}}^{o})\cong0\;,\;for\;all\;V,W\in
\Gamma(S),\;\;\Leftrightarrow\;\;\Sigma\;\;is\;III-flat\;]$
\end{proposition}

\begin{proof}
We apply Theorem \ref{theor:covar.curvat.extends}a to all $m$ terms in
(\ref{ricci}) for each case.

(a) The first $m-1$ terms have good limits on $\Sigma$, whereas $\left\langle
R({\mathcal{N}^{\prime}}^{o},\mathcal{R}^{o}){\mathcal{N}^{\prime}}%
^{o},\mathcal{R}^{o}\right\rangle $ diverges like ${\tau^{o}}^{-1}$. Thus
$Ric({\mathcal{N}^{\prime}}^{o},{\mathcal{N}^{\prime}}^{o})$ diverges like
${\tau^{o}}^{-2}$.

(b) As in (a), the first $m-1$ terms have good limits on $\Sigma$. Using
(\ref{upsilon}), we explicitly compute%

\[
\left\langle R({\mathcal{N}^{\prime}}^{o},\mathcal{R}^{o}){\mathcal{V}%
^{\prime}}^{o},\mathcal{R}^{o}\right\rangle =k^{o}\tau^{o}-\square
_{\lbrack{\mathcal{V}^{\prime}}^{o},\mathcal{R}^{o}]}{\mathcal{N}^{\prime}%
}^{o}(\mathcal{R}^{o})-\square_{\mathcal{R}^{o}}{\mathcal{N}^{\prime}}%
^{o}(\mathcal{R}^{o})({\tau^{o}}^{-1}\left\langle \mathcal{R}^{o}%
,\nabla_{{\mathcal{V}^{\prime}}^{o}}\mathcal{R}^{o}\right\rangle )
\]

\noindent($k\in C^{\infty}(M)$), and the result follows from (\ref{HII}),
(\ref{ro}) and (\ref{ro,nu-nabla}). Therefore, the possible divergence of
$Ric({\mathcal{N}^{\prime}}^{o},{\mathcal{V}^{\prime}}^{o})$ when approaching
$\Sigma$ is like ${\tau^{o}}^{-1}$.

(c) The first $m-1$ terms have good limits on $\Sigma$ (for $i=2,\ldots,m-1,$
because of $II$-flatness). Using (\ref{upsilon}), we explicitly compute%

\[
\left.
\begin{array}
[c]{l}%
\left\langle R({\mathcal{N}^{\prime}}^{o},\mathcal{R}^{o}){\mathcal{R}%
^{\prime}}^{o},\mathcal{R}^{o}\right\rangle =k^{o}\tau^{o}+\square
_{{\mathcal{R}^{\prime}}^{o}}{\mathcal{N}^{\prime}}^{o}(\mathcal{R}^{o}%
)({\tau^{o}}^{-1}\left\langle \mathcal{R}^{o},\nabla_{\mathcal{R}^{o}%
}\mathcal{R}^{o}\right\rangle )-\\
\;\;\;\;\;\;\;\;\;\;\;\;\;\;\;\;\;\;\;\;\;\;\;\;\;\;\;\;\;\;\;\;\;\;\;\;\;\;\;\;-\square
_{\mathcal{R}^{o}}{\mathcal{N}^{\prime}}^{o}(\mathcal{R}^{o})({\tau^{o}}%
^{-1}\left\langle \mathcal{R}^{o},\nabla_{{\mathcal{R}^{\prime}}^{o}%
}\mathcal{R}^{o}\right\rangle )
\end{array}
\right.
\]

\noindent($k\in C^{\infty}(M)$), and the result follows from (\ref{ro}) and
(\ref{covar.deriv.restr.depends,bis}). As in (b), the possible divergence of
$\;Ric({\mathcal{N}^{\prime}}^{o},{\mathcal{R}^{\prime}}^{o})$ when
approaching $\Sigma$ is like ${\tau^{o}}^{-1}$.

(d) The first $m-1$ terms have good limits on $\Sigma$. Using (\ref{upsilon}),
we explicitly compute $\;\left\langle R({\mathcal{R}^{\prime}}^{o}%
,\mathcal{R}^{o}){\mathcal{V}^{\prime}}^{o},\mathcal{R}^{o}\right\rangle
=k^{o}\tau^{o}$ ($k\in C^{\infty}(M)$), and the result follows.

(e) The first term diverges like ${\tau^{o}}^{-1}$, whereas the next $m-2$
terms have good limits on $\Sigma$. Moreover, $\;\left\langle R({\mathcal{R}%
^{\prime}}^{o},\mathcal{R}^{o}){\mathcal{R}^{\prime}}^{o},\mathcal{R}%
^{o}\right\rangle ={\tau^{o}}^{2}\left\langle R(\bar{\mathcal{R}}%
^{o},\mathcal{R}^{o})\bar{\mathcal{R}}^{o},\mathcal{R}^{o}\right\rangle $.
Thus $\;Ric({\mathcal{R}^{\prime}}^{o},{\mathcal{R}^{\prime}}^{o})$ diverges
like ${\tau^{o}}^{-1}$ when approaching $\Sigma\;$.

(f) The first $m-1$ terms have good limits on $\Sigma$ (for $i=2,\ldots,m-1,$
because of $II$-flatness). Using (\ref{upsilon}) and $II$-flatness, we
explicitly compute%

\[
\left\langle R({\mathcal{V}^{\prime}}^{o},\mathcal{R}^{o}){\mathcal{W}%
^{\prime}}^{o},\mathcal{R}^{o}\right\rangle =k^{o}\tau^{o}-\square
_{{\mathcal{V}^{\prime}}^{o}}{\mathcal{W}^{\prime}}^{o}(\mathcal{E}_{1}%
^{o})\square_{\mathcal{R}^{o}}\mathcal{E}_{1}^{o}(\mathcal{R}^{o})
\]

\noindent($k\in C^{\infty}(M)$), and the result follows from the fact that
$III$-flatness is equivalent to $II$-flatnes plus $\mathcal{H}$-flatness.
Thus, the possible divergence of $\;Ric({\mathcal{V}^{\prime}}^{o}%
,{\mathcal{W}^{\prime}}^{o})$ when approaching $\Sigma$ is like ${\tau^{o}%
}^{-1}\;\;\;\;\;\;$
\end{proof}

\begin{remark}
The results in Proposition \ref{prop:ricci curvat.extends} sharply contrasts
with what happens in the case of transverse radical; for then, $III$-flatness
is equivalent (\cite{kossow2}, Theorem 3c) to the fact that $\;Ric(\mathcal{A}%
^{o},\mathcal{B}^{o})\cong0$, for all $\mathcal{A},\mathcal{B}\in
\frak{X}(M)\;\;\;\;\;\endproof$
\end{remark}

\begin{proposition}
\label{prop:ricci curvat.restricts} Let $(M,g)$ be a transverse
Riemann-Lorentz space with tangent radical and singular hypersurface $\Sigma$.
Let be $A,B\in\frak{X}_{\Sigma}\,$. Then, even in the $III$-flat case,
$Ric(A,B)$ is not well-defined.
\end{proposition}

\begin{proof}
We apply Theorem \ref{theor:covar.curvat.restricts} to all $m$ terms in
(\ref{ricci}) for the cases (d) and (f) in Proposition \ref{prop:ricci
curvat.extends} (the only cases we need to check).

Case (d): $\left\langle R(R,N)V,N\right\rangle $ is never well-defined,
whereas $\left\langle R(R,E_{\lambda})V,E_{\lambda}\right\rangle $ is
well-defined if $\Sigma\,$ is $II$-flat ($\lambda=2,\ldots,m-1$) and
$\left\langle R(R,R)V,R\right\rangle $ is always well-defined (actually, it
vanishes). Thus $Ric(R,V)$ is not well-defined.

Case (f): $\left\langle R(V,N)W,N\right\rangle $ is always well-defined and
$\left\langle R(V,E_{\lambda})W,E_{\lambda}\right\rangle $ is well-defined if
$\Sigma\,$ is $II$-flat ($\lambda=2,\ldots,m-1$). However, although
\linebreak ${\tau^{o}}^{-1}\left\langle R({\mathcal{V}^{\prime}}%
^{o},\mathcal{R}^{o}){\mathcal{W}^{\prime}}^{o},\mathcal{R}^{o}\right\rangle
\cong0\;$ if $\Sigma\,$ is $III$-flat, computation using (\ref{upsilon}) shows
that $\;({\tau^{o}}^{-1}\left\langle R({\mathcal{V}^{\prime}}^{o}%
,\mathcal{R}^{o}){\mathcal{W}^{\prime}}^{o},\mathcal{R}^{o}\right\rangle
)\mid_{\Sigma}\;$ depends on the extensions $\mathcal{V}^{\prime}%
,\mathcal{W}^{\prime}$ of $V,W$. Thus, $Ric(V,W)$ is not well-defined $\;\;\;\;\;$
\end{proof}

\bigskip

\subsection{The tangential connection}

When the singular hypersurface $\Sigma$ is $II$-flat, Proposition
\ref{prop:covar.deriv.}b induces the following natural connection on $\Sigma$
(we denote it by $\;\nabla^{\Sigma}\,$ and we call it the \emph{tangential
connection})%

\[
\nabla^{\Sigma}:\frak{X}(\Sigma)\times\frak{X}(\Sigma)\rightarrow
\frak{X}(\Sigma),(X,Y)\mapsto\;\nabla_{X}^{\Sigma}Y:=\nabla_{\mathcal{X}%
}\mathcal{Y}\mid_{\Sigma}-\nu(\nabla_{\mathcal{X}}\mathcal{Y}\mid_{\Sigma
})N\;,
\]

\noindent where $\mathcal{X},\mathcal{Y} \in\frak{X}(M)$ are the
\emph{canonical} extensions of $X,Y$. Using (\ref{covar.deriv.}), we obtain:%

\begin{equation}
\nabla^{\Sigma}_{X} Y = \Sigma^{m-1}_{\lambda=2} \square_{X} Y (E_{\lambda})
E_{\lambda}+ (\tau^{-1} \left\langle \nabla_{\mathcal{X}} \mathcal{Y}%
,\mathcal{R} \right\rangle )\mid_{\Sigma}\; R \;\;.\; \label{tangential1}%
\end{equation}

As a consequence of (\ref{ro,nu-nabla}), \emph{the main radical vectorfield
$R$ becomes $\nabla^{\Sigma}$-geodesic}.

\bigskip

Equations (\ref{ro}), (\ref{Xtau}) and (\ref{nabla antysim. restr.}) lead to:%

\[
\left.
\begin{array}
[c]{l}%
X(\rho(Y)) = X(\tau^{-1} \left\langle \mathcal{R},\mathcal{Y}\right\rangle ) =
(\tau^{-1} \left\langle \nabla_{\mathcal{X}} \mathcal{R},\mathcal{Y}
\right\rangle )\mid_{\Sigma}\;+ (\tau^{-1} \left\langle \mathcal{R}%
,\nabla_{\mathcal{X}} \mathcal{Y} \right\rangle )\mid_{\Sigma}\;,\;\Rightarrow
\\
\\
\Rightarrow\;(\tau^{-1} \left\langle \nabla_{\mathcal{X}} \mathcal{R}%
,\mathcal{Y} \right\rangle _{Ant})\mid_{\Sigma}\;=\\
\\
\;\;\;\;\;\;\;\;\;\; = \frac{1}{2} \{X(\rho(Y)) - Y(\rho(X))\} - (\tau^{-1}
\left\langle \mathcal{R},\nabla_{\mathcal{X}} \mathcal{Y} \right\rangle
_{Ant})\mid_{\Sigma}\;= \frac{1}{2} d\rho(X,Y)\;,
\end{array}
\right.
\]

\smallskip

\noindent from which we obtain (using (\ref{tangential1}), (\ref{main}) and
(\ref{main admissible})):%

\[
\left.
\begin{array}
[c]{c}%
\nabla^{\Sigma}_{X} Y = D^{S}_{X} Y + \{ X(\rho(Y)) - (\tau^{-1} \left\langle
\nabla_{\mathcal{X}} \mathcal{R},\mathcal{Y} \right\rangle )\mid_{\Sigma}\}\;
R =\\
\\
= \tilde{D}_{X} Y - (\tau^{-1} \left\langle \nabla_{\mathcal{X}}
\mathcal{R},\mathcal{Y} \right\rangle )\mid_{\Sigma}\; R = \tilde{D}_{X} Y -
\frac{1}{2} d \rho(X,Y) R - (\tau^{-1} \left\langle \nabla_{\mathcal{X}}
\mathcal{R},\mathcal{Y} \right\rangle _{Sim})\mid_{\Sigma}\; R =\\
\\
= \dot{D}_{X} Y - (\tau^{-1} \left\langle \nabla_{\mathcal{X}} \mathcal{R}%
,\mathcal{Y} \right\rangle _{Sim})\mid_{\Sigma}\; R \;\;.
\end{array}
\right.
\]

It follows from Theorem \ref{theor:admissible} that \emph{$\nabla^{\Sigma}$ is
an admissible, metric connection} and that \emph{all admissible connections on
$\Sigma$ have the same covariant curvature as $\nabla^{\Sigma}$}. We finally
compute that covariant curvature:

\begin{theorem}
\label{theor:covar.curvat.tangential} Let $(M,g)$ be a transverse
Riemann-Lorentz space with tangent radical and singular, $II$-flat
hypersurface $\Sigma$. Let $R^{\Sigma}$ be the curvature of the tangential
connection $\nabla^{\Sigma}$ and let be $X,Y,Z,T\in\frak{X}(\Sigma)$. Then it
holds (Gauss formula):%

\[
\left\langle R^{\Sigma}(X,Y)Z,T\right\rangle =\left\langle
R(X,Y)Z,T\right\rangle -\det\left(
\begin{array}
[c]{c}%
\mathcal{H}(X,Z)\;\;\mathcal{H}(Y,Z)\\
\mathcal{H}(X,T)\;\;\mathcal{H}(Y,T)
\end{array}
\right)  \;\;.
\]
\end{theorem}

\textbf{Proof} By $II$-flatness and Theorem \ref{theor:covar.curvat.restricts}%
, $\left\langle R(X,Y)Z,T\right\rangle $ is well-defined. We compute it using
the canonical extensions $\mathcal{X},\mathcal{Y},\mathcal{Z}$ of $X,Y,Z$.
Taking into account that (by (\ref{ro,nu-nabla})):%

\[
\left.
\begin{array}
[c]{c}%
\nabla_{\mathcal{Y}} \mathcal{Z}\mid_{\Sigma}\; = \nabla^{\Sigma}_{Y} Z +
\mathcal{H} (Y,Z)N\;,\;\;\Rightarrow\\
\\
\Rightarrow\;\;\square_{X} (\nabla_{\mathcal{Y}} \mathcal{Z}\mid_{\Sigma
}\;)(T) = \square_{X} (\nabla^{\Sigma}_{Y} Z)(T) - \mathcal{H} (Y,Z)
\mathcal{H}(X,T) =\\
\\
= \left\langle \nabla^{\Sigma}_{X} (\nabla^{\Sigma}_{Y} Z), T \right\rangle -
\mathcal{H} (Y,Z) \mathcal{H}(X,T)\;\;,
\end{array}
\right.
\]

\noindent we obtain:%

\[
\left.
\begin{array}
[c]{l}%
\left\langle R(X,Y)Z,T \right\rangle = \square_{X} (\nabla_{\mathcal{Y}}
\mathcal{Z}\mid_{\Sigma}\;)(T) - \square_{Y} (\nabla_{\mathcal{X}}
\mathcal{Z}\mid_{\Sigma}\;)(T) - \square_{[X,Y]}Z (T) =\\
\\
\;\;\;\;\; = \left\langle \nabla^{\Sigma}_{X} (\nabla^{\Sigma}_{Y} Z), T
\right\rangle - \mathcal{H} (Y,Z) \mathcal{H}(X,T) -\\
\\
\;\;\;\;\;\;\;\;\;\;\;\;\;\;\;\;\;\;\;\; - \left\langle \nabla^{\Sigma}_{Y}
(\nabla^{\Sigma}_{X} Z), T \right\rangle + \mathcal{H} (X,Z) \mathcal{H}(Y,T)
- \left\langle \nabla^{\Sigma}_{[X,Y]} Z,T \right\rangle =\\
\\
\;\;\;\;\; = \left\langle R^{\Sigma}(X,Y)Z,T \right\rangle + \det\left(
\begin{array}
[c]{c}%
\mathcal{H}(X,Z) \;\; \mathcal{H}(Y,Z)\\
\mathcal{H}(X,T) \;\; \mathcal{H}(Y,T)
\end{array}
\right)  \;\;\;\;\;\;\;\;\;\;\endproof
\end{array}
\right.
\]

\bigskip

\bibliographystyle{plain}
\bibliography{radical1}
\end{document}